\documentclass[a4paper, 12pt]{article}

\usepackage[sort&compress]{natbib}
\bibpunct{(}{)}{;}{a}{}{,} 

\usepackage{amsthm, amsmath, amssymb,   mathrsfs, multirow, url, subfigure}
\usepackage{graphicx} 
\usepackage{ifthen} 
\usepackage{amsfonts}
\usepackage[usenames]{color}
\usepackage{fullpage}
\usepackage{enumitem}

\usepackage{mathtools}

\theoremstyle{plain}

\theoremstyle{definition}

\theoremstyle{remark}

\newtheorem{cond}{Condition}

\newcommand{\prob}{\mathsf{P}} 
\newcommand{\E}{\mathsf{E}}
\newcommand{\V}{\mathsf{V}}

\newcommand{\nm}{{\sf N}}

\newcommand{\icomp}{\mathrm{i}}

\renewcommand{\S}{\mathcal{S}}

\renewcommand{\phi}{\varphi} 
\newcommand{\eps}{\varepsilon}

\newcommand{\event}{\mathscr{E}}

\newcommand{\iid}{\overset{\text{\tiny iid}}{\,\sim\,}}
\newcommand{\ind}{\overset{\text{\tiny ind}}{\,\sim\,}}

\theoremstyle{plain} 
\newtheorem{thm}{Theorem}

\newtheorem{lem}{Lemma}

\newcommand{\ceiling}[1]{\lceil #1 \rceil}

\newcommand{\K}{\mathcal{K}}

\title{Gibbs posterior inference on a L\'evy density under discrete sampling}
\author{
Zhe Wang\footnote{Department of Statistics, North Carolina State University; {\tt zwang54@ncsu.edu}, {\tt rgmarti3@ncsu.edu}} \quad and \quad Ryan Martin$^*$
}
\date{\today}

\begin{document}

\maketitle

\begin{abstract}
In mathematical finance, L\'evy processes are widely used for their ability to model both continuous variation and abrupt, discontinuous jumps.  These jumps are practically relevant, so reliable inference on the feature that controls jump frequencies and magnitudes, namely, the {\em L\'evy density}, is of critical importance.  A specific obstacle to carrying out model-based (e.g., Bayesian) inference in such problems is that, for general L\'evy processes, the likelihood is intractable.  To overcome this obstacle, here we adopt a {\em Gibbs posterior} framework that updates a prior distribution using a suitable loss function instead of a likelihood.  We establish asymptotic posterior concentration rates for the proposed Gibbs posterior.  In particular, in the most interesting and practically relevant case, we give conditions under which the Gibbs posterior concentrates at (nearly) the minimax optimal rate, adaptive to the unknown smoothness of the true L\'evy density. 

\smallskip

\end{abstract}

\section{Introduction}
\label{S:intro}

Mathematical finance is largely focused on modeling the fluctuations in asset prices over time.  Classical models, such as that of \citet{MR3363443}, assume that the sample paths of an asset are continuous. Recent technological advances, however, have made high-frequency trading possible, which means prices can now change drastically and effectively instantly. Therefore, {\em L\'evy processes} \citep[e.g.,][]{MR2512800, MR1406564, MR3185174, MR2042661, obn.levy.collection.2001} 
have become an essential building block for models that allow for flexible jump behavior in the sample paths.
More specifically, a (real-valued) L\'evy process $X = \{X(t): t \geq 0\}$ is a collection of random variables, where $t$ is often interpreted as time, with the three key properties: independent and stationary increments, right-continuous with left-hand limits (i.e., c\`adl\`ag), and no fixed jump times; see Section~\ref{SS:model} and the references given there.  Important special cases include Brownian motion with drift and compound Poisson processes, corresponding to the ``no jumps'' and ``only jumps'' extremes, respectively.  By considering more general L\'evy processes,  one's model can simultaneously accommodate both jump discontinuities and asymmetric/heavy-tailed distributions for $X(t+s) - X(t)$.  

The general L\'evy process is characterized by a triplet $(\mu, \sigma, \Psi)$, where $\mu$ and $\sigma$ are coefficients related to the continuous part of the sample path (basically a Brownian motion with drift), and $\Psi$ is a sigma-finite measure on $\mathbb{R}\setminus\{0\}$---the {\em L\'evy measure}---that determines the frequency and magnitudes of the jump discontinuities.  Roughy, for each (measurable) subset $A \subset \mathbb{R}\setminus\{0\}$, the L\'evy measure $\Psi(A)$ represents the average number of jumps per unit time whose magnitudes fall into $A$.  Since the jumps would be relevant to so-called ``black swan'' events \citep[e.g.,][]{taleb2007black}, which are crucial for risk assessments, our focus is estimation of/inference on $\Psi$.  More specifically, as it is common in applications to assume that $\Psi$ has a density with respect to Lebesgue measure, our object of interest is the {\em L\'evy density} $\psi$ that satisfies $\Psi(A) = \int_A \psi(x) \, dx$. 


For statistical inference on $\Psi$ or $\psi$, there are two key challenges.  First, only the jump frequencies and sizes in the $X$ sample path are directly related to the L\'evy measure, but these are only observable under continuous-time sampling, which is rarely possible in practice.  Instead, typically $X$ can be observed only at discrete time points, $t_0 < t_1 < \cdots < t_n$, which means that jump features are not directly observable.  That is, $X(t_i) - X(t_{i-1})$ is almost surely non-zero, but it is impossible to separate the part corresponding to jumps in $[t_{i-1},t_i]$ from that corresponding to continuous variation therein.  However, it is possible, at least asymptotically, to recover the relevant information needed to infer the L\'evy measure under a {\em high-frequency} sampling scheme \citep{sahalia.jacod.book}, as we adopt here.  More specifically, we assume that the length $t_n$ of the observation window and the spacing, denoted by $\Delta=\Delta_n$, are approaching $\infty$ and $0$, respectively, as $n \to \infty$. The intuition behind such a sampling scheme is that we nearly observe the jumps when the sampling frequency is high, and eventually there will be a large number of such jumps when the time horizon is large, hence $\Psi$ can be recovered.  

The second challenge is that the general L\'evy process does not determine a closed-form likelihood function for $\Psi$ as a function of the observable $\{X(t_0), \ldots, X(t_n)\}$ that can be used in model-based inference.  Indeed, the distribution of $X(t)$ is directly connected to the parameters $(\mu, \sigma, \Psi)$ through its characteristic function and the L\'evy--Khintchine formula, not through a density function.  Consequently, at least in the general case, there is no convenient likelihood function in $\psi$ available to carry out maximum likelihood or Bayesian inference.  For this reason, M-estimation techniques, as discussed in Section~\ref{S:background}, are most commonly used for estimation of and inference on the L\'evy density.

The lack of an explicit likelihood motivates us to consider a model-free Bayesian-like approach, the so-called {\em Gibbs posterior}, as advocated for in \citet{MR2291497, MR2241190, MR3557191}. Suppose the quantity of interest minimizes some appropriately defined risk function.  Then the Gibbs posterior treats a scaled empirical risk function as a negative log-likelihood and applies Bayes's formula to obtain a posterior distribution.  For the present context, it has been shown \citep[e.g.,][]{figueroa2009nonparametric} that the L\'evy density 
can be viewed as the minimizer of a limiting risk or contrast function, which is what inspires the aforementioned M-estimation techniques for inference.  So, the Gibbs posterior framework for estimation of and uncertainty quantification for $\psi$ seems quite natural.  This is not straightforward, however, since $\psi$ is an infinite-dimensional quantity.  


The present paper's focus is on the construction of a suitable Gibbs posterior distribution that can be used for probabilistic inference on $\psi$.  Note, again, that a proper Bayesian approach for inference on $\psi$ is out of reach, at least in the general case, because there is no likelihood function available.  Theoretical support for the proposed Gibbs posterior comes in the form of a demonstration that it achieves the optimal, asymptotic concentration rate properties in various settings, under appropriate conditions on the true L\'evy density, the spacing $\Delta$, etc.  In particular, the highlight of the paper (Theorem~\ref{thm:K.unknown}) establishes that, if the true L\'evy density, $\psi^\star$, is order-$\alpha$ smooth with $\alpha$ unknown, then the proposed Gibbs posterior concentrates, asymptotically and adaptively, around $\psi^\star$ at (nearly) the $\alpha$-dependent minimax optimal rate {\em for estimators with access to the full sample path $X$}.  That is, remarkably, there is virtually no loss of efficiency for not knowing the smoothness of $\psi^\star$ or for not having access to the full sample path $X$.  

The remainder of the paper is organized as follows.  In Section~\ref{S:background}, we review the relevant details of the general L\'evy process model, as well as existing approaches for estimating the corresponding L\'evy density.  Section~\ref{S:gibbs} presents our proposed approach, wherein a suitable loss (or contrast) function is used in place of a likelihood to update a prior distribution, making ``model-free'' posterior inference on $\psi$ possible.  In the most general case, we introduce a {\em complexity index} and a corresponding hierarchical prior, in the spirit of \citet{MR3091697} and \citet{MR3426318}, which makes it possible for the posterior to adapt to the unknown smoothness of the underlying L\'evy density.  The main results are in Section~\ref{S:theory}, where we present three asymptotic theorems, in increasing order of generality, culminating in Theorem~\ref{thm:K.unknown} which establishes the Gibbs posterior's adaptive and (nearly) minimax optimal concentration rate.  We also show, in Theorem~\ref{thm:no.overfit}, that the marginal Gibbs posterior for the complexity index will concentrate on a range of values compatible with the oracle complexity that depends on the smoothness of $\psi^\star$.  
A remarkable feature of the proposed Gibbs formulation is that, despite the complexities of the L\'evy process and the infinite-dimensionality of the L\'evy density, posterior computation is straightforward.  In Section~\ref{S:numerical} we present first a simple Monte Carlo algorithm for simulating the Gibbs posterior for $\psi$ and then give a brief numerical illustration.  Concluding remarks are given in Section~\ref{S:discuss} and proofs of all the technical results can be found in the Appendices.

\section{Background}
\label{S:background}

\subsection{Model}
\label{SS:model}

Let $X = \{X(t): t \geq 0\}$ be a real-valued L\'evy process  and write $\prob$ for its probability distribution or law. That is, $X$ satisfies the following properties: 
\begin{itemize}
\item {\em fixed starting value}, i.e., without loss of generality, $X(0)=0$  almost surely;
\item {\em independent increments}, i.e., for any pair of non-overlapping time intervals, say, $[s_1,t_1]$ and $[s_2,t_2]$, the random variables $X(t_1) - X(s_1)$ and $X(t_2) - X(s_2)$ are independent; 
\item {\em stationary increments}, i.e., for any $s<t$, the random variables $X(t)-X(s)$ and $X(t-s)$ have the same distribution; 
\item and {\em stochastic continuity}, i.e., $\prob\{ |X(t+h) - X(t)| > \eps) \to 0$ as $h \to 0$, for any $\eps > 0$ and $t>0$.
\end{itemize}
The two most common examples of L\'evy processes are Brownian motion, which has continuous sample paths almost surely, and compound Poisson processes, which have a finite number of jumps on every bounded time interval almost surely. Of course, there are other examples, including ones with infinite jump activity \citep[e.g.,][]{madan1990variance}.  In general, the distribution $\prob$ of a general L\'evy process $X$ is determined by the three parameters, $(\mu, \sigma, \Psi)$, as elegantly described via the L\'evy--Khintchine representation  \citep{levy1934integrales} of the characteristic function of $X(t)$, for $t \geq 0$, i.e., 
\begin{equation}
\label{eq:LK}
\E \{e^{\icomp \gamma X(t)}\} =\exp\Bigl[ t\Bigl\{ \mu \icomp \gamma-\frac{\sigma^2\gamma^2}{2}+\int(e^{\icomp \gamma x}-1-\icomp \gamma x 1_{|x|<1}) \, \Psi(dx)\Bigr\} \Bigr], \quad \gamma \in \mathbb{R},
\end{equation}
where $\icomp=(-1)^{1/2}$ denotes the complex root, $1_A$ denotes the indicator function of an event $A$, and $\E$ denotes expectation with respect to the distribution $\prob$ of $X$.  Note that, while the characteristic function fully determines the underlying distribution, it is not generally possible to work invert the characteristic function/Fourier transform and obtain a closed-form expression for the density and, in turn, a likelihood function, depending explicitly on the to-be-estimated $\Psi$.  That this problem generally lacks a workable likelihood function makes model-based estimation and inference a challenge, as described in Section~\ref{SS:old} below.  The L\'evy--Khintchine representation above requires 
\begin{equation}
\label{eq:integrable}
\int \min(x^2, 1) \, \Psi(dx) < \infty.
\end{equation}
That \eqref{eq:integrable} holds for the true L\'evy measure being estimated is a standing assumption made throughout the paper.  

A more descriptive characterization of the L\'evy process structure is given by the celebrated L\'evy--It\^o decomposition \citep{MR14629}.  The first step is to split the process into its continuous and jump parts as 
\[ X(t) = C(t) + J(t). \]
The continuous part, $C(\cdot)$, can be shown to be a Brownian motion with a linear drift of slope $\mu$ and volatility coefficient $\sigma$.  Naturally, the jump part, $J(\cdot)$, is more complicated. It\^o  showed, first, that the jump part can be written as 
\[ J(t) = \sum_{s: 0 \leq s \leq t} \{X(s) - X(s-)\}. \]
and, second, that it can be further decomposed into a the sum of a compound Poisson process and a certain limit of a compensated compound Poisson process.  Start by defining the Poisson random measure 
\[ N(t, A) = |\{ s \in [0,t]: X(s) - X(s-) \in A\}|, \]
where $|S|$ denotes the cardinality of a finite set $S$.  The idea is that, for a fixed $t$, $A \mapsto N(t,A)$ is a random measure and, for fixed $A$, $t \mapsto N(t,A)$ is a Poisson process with rate $\Psi(A)$.  Then, for any $h > 0$, the jump part of the L\'evy process can be written as 
\[ J(t) = \int_{|x| > h} x \, N(t,dx) + \int_{|x| < h} x \, \{N(t,dx) - t \Psi(dx)\}, \]
where the latter is defined as $\lim_{\eps \to 0} \int_{\eps < |x| < h} x \, \{N(t,dx) - t \Psi(dx)\}$, in mean square. This expresses $J$ as a an accumulation of finitely many ``large'' jumps, as determined by the compound Poisson process, and potentially infinitely many ``small'' jumps, as determined by the compensated compound Poisson process.  
From this ``$X=C+J$'' decomposition, it is clear that Brownian motion and compound Poisson processes are important special cases of L\'evy processes.  In the latter case, if $\lambda$ is the jump intensity of the Poisson process and $G$ is the jump size distribution, then the L\'evy measure is given by $\Psi(dx) = \lambda G(dx)$.  Note that $\Psi$ is generally not a probability measure. 

A distinguishing feature of L\'evy processes, and what makes them attractive models for use in applications, is the jumps.  Since the features of these jumps are controlled by $\Psi$, which is generally unknown in practice, it would be desirable to estimate or make inference on $\Psi$ (or functionals thereof) based on data. If the entire sample path, $t \mapsto X(t)$, were observable, then the jump frequencies and sizes would be observable as well and these observations could be used directly to estimate features of $\Psi$ relatively easily.  In applications, however, it is rare for the entire sample path to be accessible; a more realistic situation is where one observes the sample path a some large---but ultimately finite---set of time points, say, $0 = t_0 < t_1 < \cdots < t_n$.  This means that the jump sizes and frequencies, the quantities that contain the relevant information about $\Psi$, are not directly observable.  In such cases, only in an asymptotic regime under high-frequency sampling can the features of $\Psi$ be estimated accurately.  See Section~\ref{SS:old} below for more specifics about the high-frequency sampling scheme adopted here.

\subsection{Existing approaches}
\label{SS:old}

As mentioned above, the complexity of a general L\'evy process makes model-based estimation and inference on $\Psi$ a challenge.  For that reason, model-based approaches---Bayesian in particular---have focused on special cases of the general L\'evy process, such as compound Poisson processes \citep{MR3356921, MR3769832, MR4013745}, where a likelihood function readily is available. More recently, \citet{belomestny2018nonparametric} developed a nonparametric Bayes approach for inference in another special case, namely, gamma-type L\'evy subordinators, in which they use a data-augmentation scheme to avoid the intractable likelihood function.  Unfortunately, to our knowledge, their approach cannot be readily extended to the general case under consideration here.

Without a likelihood, how can one estimate and make inference on $\Psi$?  As indicated above, a first requirement is that, for general L\'evy processes, we need a high-frequency sampling scheme in order to separate the information relevant to the jump part of the process from that relevant to the continuous part.  Recall that the process $X$ is to be observed at the times points $0=t_0 < t_1 < \cdots < t_n$, a collection indexed by $n$.  Following \citet{MR2565560}, \citet{figueroa2009nonparametric}, \citet{adaptive.levy.2011}, and others, we consider a high-frequency sampling scheme where $t_n \to \infty$ and $\Delta_n = \max_{i=1,\ldots,n} (t_i-t_{i-1}) \to 0$, simultaneously, as $n \to \infty$.  In contrast, a {\em low-frequency} sampling scheme corresponds to keeping $\Delta_n$ fixed while $t_n \to \infty$.  

For simplicity, and with virtually no loss of generality, we will assume throughout the paper that the observation times are equally spaced in the interval $[0,t_n]$, so that $t_i = i \Delta$, for $i=0,1,\ldots,n$. Then the increments $Y_i = X(t_i) - X(t_{i-1})$, for $i=1,\ldots,n$ are independent and identically distributed (iid).  By stationarity, $X_\Delta = X(\Delta)$ has the same distribution as the $Y_i$'s.  A critical observation made by \citet{figueroa2009nonparametric} is that, for any bounded and continuous function $f$, vanishing in a neighborhood of the origin, the following small-time result holds:
\begin{equation}
\label{eq:jose.smalltime}
\Delta^{-1} \, \E f(X_\Delta) = \int f(x) \, \psi^\star(x) \, dx + O(\Delta), \quad \Delta \to 0. 
\end{equation}
This motivates a powerful and flexible {\em projection estimator}, inspired by the method of sieves \citep[e.g.,][]{MR599175} and the model selection literature \citep[e.g.,][]{MR1462939}. That is, for a positive integer $K \geq 1$, consider a $K$-dimensional linear space $\S_K \subset L_2(D)$ of functions that are bounded, continuous, and vanishing outside the interval $D$, which itself is bounded away from the origin.  More specifically, let 
\[ \S_K = \text{span}\{f_{Kk}: k=1,\ldots,K\}, \]
where the $f_{Kk}$'s are bounded and continuous functions, vanishing outside of $D$.  Without loss of generality, let $f_{K1},\ldots,f_{KK}$ be orthonormal and form a basis for $\S_K$.  Then the $L_2$ projection of the L\'evy density $\psi^\star$ onto $\S_K$ is given by 
\[ \psi_K^\perp(x) = \sum_{k=1}^K \theta_{Kk}^\perp \, f_{Kk}(x), \]
where 
\[ \theta_{Kk}^\perp = \int_D f_{Kk}(x) \, \psi^\star(x) \, dx, \quad k=1,\ldots,K. \]
This, together with \eqref{eq:jose.smalltime}, suggests a projection-based estimator, namely, 
\[ \hat\psi_K(x) = \sum_{k=1}^K \hat\theta_{Kk} \, f_{Kk}(x), \]
where 
\[ \hat\theta_{Kk} = \frac{1}{t_n} \sum_{i=1}^n f_{Kk}(Y_i), \quad k=1,\ldots,K. \]
It is intuitively clear from \eqref{eq:jose.smalltime} that, at least for fixed $K$, the projection estimator $\hat\psi_K$ converges $\psi_K^\perp$ as $n \to \infty$, provided that $\Delta \to 0$ and $t_n \to \infty$.  This intuition can be made rigorous, e.g., Proposition~3.4 of \citet{figueroa2009nonparametric} states that, if the $f_{Kk}$'s have derivatives uniformly bounded in $k$  
\[ \E \| \hat\psi_K - \psi_K^\perp \|_{L_2(D)}^2 \leq c_K \, K \, t_n^{-1}, \]
where $\|\cdot\|_{L_2(D)}$ denotes the usual $L_2(D)$-norm, and $c_K$ is a constant that depends on $\psi^\star$ and on $K$ through certain features of the $f_{Kk}$'s.  Therefore, if $K$ is fixed, then $\hat\psi_K$ converges to $\psi_K^\perp$ in $L_2(D)$ at rate $t_n^{-1}$.  

However, the objective is not to estimate $\psi_K^\perp$, it is to estimate $\psi^\star$.  This depends on how well $\psi^\star$ can be approximated by functions in $\S_K$ which, in turn, depends on how large $K$ is.  By the Pythagorean identity, 
\[ \E\|\hat\psi_K - \psi^\star\|_{L_2(D)}^2 = \E\|\hat\psi_K - \psi_K^\perp\|_{L_2(D)}^2 + \|\psi_K^\perp - \psi^\star\|_{L_2(D)}^2. \]
The first term on the right-hand side, what \citet{figueroa2009nonparametric} calls the ``variance term,'' has an upper bound as a function of $K$ and $t_n$ as above.  The second term, which represents the bias of the model $\S_K$, or the approximation error, also has known upper bounds, depending on certain features of $\S_K$ and the smoothness of $\psi^\star$.  Suitably balancing the bounds on the two terms on the right-hand side above establishes the convergence rate of the projection estimator.  To make this approach practically useful, one needs a data-driven solution to the unknown-$K$ problem.  A natural idea is to somehow produce an estimator ``$\widehat K$'' of the unknown complexity, and return $\hat\psi = \hat\psi_{\widehat K}$ as the final, plug-in projection estimator.  See \citet{adaptive.levy.2011} for a specific choice of $\widehat K$ and convergence properties of the corresponding plug-in estimator $\hat\psi$.

\section{Gibbs posterior distribution}
\label{S:gibbs}


The classical Bayesian formulation is designed for cases in which the data and quantities of interest are linked through a likelihood function, in particular, when the quantities of interest are parameters of a posited statistical model.  In other cases, however, the quantities of interest are not naturally interpreted as model parameters or, more generally, there may not be a likelihood function that links the data to the quantities of interest.  A more general formulation is as described in \cite{MR3082220}, where the quantity of interest is connected to the data through a loss function or, more specifically, where the quantity of interest is defined as the minimizer of a risk/expected loss function. In such cases, a so-called {\em Gibbs posterior} formulation is recommended; see, e.g., \cite{MR3557191, MR4095335, syring2021gibbs}.  In this framework, (a multiple of) an empirical version of that risk function is treated like a negative log-likelihood function and then Bayes's formula is applied as usual.  As discussed briefly in Section~\ref{S:intro}, the advantage of this construction and perspective is that one can directly target quantities of interest without imposing extra model assumptions that may lead to model misspecification bias. For our present L\'evy density estimation problem, a slight variation on the typical Gibbs posterior construction is needed, as we describe next. 

A slightly different perspective on the projection estimators presented in Section~\ref{SS:old} above is as a sort of ``empirical risk minimizer.''  In \citet{MR2565560}, \citet{figueroa2009nonparametric}, and elsewhere, this ``empirical risk'' function is referred to as a {\em contrast function}.  Names aside, if we define 
\[ R_{n,K}(\theta) = -2\langle \theta, \hat\theta_K \rangle + \|\theta\|^2, \quad \theta \in \Theta_K \subseteq \mathbb{R}^K, \]
then it is immediately clear that 
\[ \hat\theta_K = \arg\min_{\theta \in \Theta_K} R_{n,K}(\theta). \]
Note that $R_{n,K}$ depends on $\Delta$, but since $\Delta$ depends on $n$ too, we try to keep the notation simple and use only a single subscript ``$n$'' to indicate the dependence on data.  Then there is a corresponding ``population/limiting risk'' which is obtained by taking expectation of $R_{n,K}$, then a limit as $\Delta \to 0$, and then applying \eqref{eq:jose.smalltime}, which gives 
\[ R_K(\theta) = -2\langle \theta, \theta_K^\perp \rangle + \|\theta\|^2, \quad \theta \in \Theta_K. \]
Then, of course, $\theta_K^\perp$ is the minimizer of this population risk.  With this perspective, we are now ready to construct the proposed Gibbs posterior for $\psi$.  

One of the advantages of the Gibbs posterior over the basic project estimators is its ability to incorporate available prior information about $\psi$. Here, since $\psi$ is being represented by elements in the linear space $\S_K$, or a union thereof, it will be convenient for us to express the prior distribution in terms of the pair $(K, \theta_K)$.  This representation suggests a hierarchical prior formulation, i.e., a marginal prior mass function $\tilde\pi$ for $K$ and a conditional prior $\tilde\Pi_K$ for $\theta_K$, given $K$.  This induces a prior $\Pi$ for $\psi$, namely, 
\[ \Pi(A) = \sum_K \tilde\pi(K) \, \tilde\Pi_K(\theta_K \in \Theta_K: \psi_{\theta_K} \in A), \quad A \subseteq L_2(D), \]
where $\psi_{\theta_K} = \sum_{k=1}^K \theta_{Kk} f_{Kk}$.  There may not be much useful prior information available to use for constructing $\tilde\Pi_K$, in which case, one can choose a relatively diffuse prior.  However, $K$ is very much related to the smoothness of $\psi^\star$---see below---and one may have some justifiable beliefs about this smoothness which can be used to inform the prior $\tilde\pi$ for $K$.

The Gibbs posterior is obtained by combining this prior information with the information in the observed data---the increments $X(t_i)-X(t_{i-1})$, $i=1,\ldots,n$---as quantified by the empirical loss functions $R_{n,K}$ as follows:
\begin{equation}
\label{eq:gibbs}
\Pi_n(A) = \frac{\sum_K \tilde\pi(K) \, \int_{\{\theta_K \in \Theta_K: \psi_{\theta_K} \in A\}} e^{-\omega t_n R_{n,K}(\theta_K)} \, \tilde\Pi_K(d\theta_K)}{\sum_K \tilde\pi(K) \int_{\Theta_K} e^{-\omega t_n R_{n,K}(\theta_K)} \, \tilde\Pi_K(d\theta_K)}, \quad A \subseteq L_2(D). 
\end{equation}
The parameter $\omega > 0$ is call the {\em learning rate} and its choice is critical to the finite-sample performance of the Gibbs posterior; see, e.g.,  \cite{MR3724979, MR3949315, MR3949316}. Here our focus is on the asymptotic properties of the Gibbs posterior, so we treat $\omega$ as a fixed constant. 

In the machine learning literature, a Gibbs posterior is often referred to as a ``randomized estimator,'' i.e., a distribution that can generate samples of the quantity of interest that tend to be close, in some sense, to the risk minimizer.  In addition, a Gibbs posterior can be used exactly like one would typically use a Bayesian posterior.  That is, the posterior mean could be used as an estimator of $\psi^\star$.  Similarly, other features of $\psi^\star$ with appropriate summaries of $\Pi_n$.  The Gibbs posterior also provides uncertainty quantification.  For example, at least in principal, probabilities for any assertions ``$\psi^\star \in A$'' can be evaluated, which could be used for testing hypotheses, etc.  Moreover, $100(1-\alpha)$\% credible sets for $\psi^\star$ can be readily obtained, which are of the form 
\[ \{\psi: d(\psi, \hat\psi) \leq c_n(\alpha)\}, \]
where $d$ is a metric, $\hat\psi$ is the Gibbs posterior mean, and $c_n(\alpha)$ is the upper-$\alpha$ quantile of the marginal Gibbs posterior for $d(\psi, \hat\psi)$, when $\psi \sim \Pi_n$.  For example, if $d$ were the supremum norm on $D$, then the above display gives a uniform credible band for $\psi^\star$.  

It turns out that the choice of how to quantify prior uncertainty about $K$ is quite important.  For example, a naive strategy would be to just take $K$ to be a fixed constant, which amounts to $\tilde\pi$ to assign probability 1 to a single value.  That simplifies the Gibbs posterior significantly, since then there would be no sums in \eqref{eq:gibbs}.  But there is a price to pay for the added simplicity in terms of the model flexibility.  For that reason, we will refer to $K$ below as the model's {\em complexity index} and how the user's choice to handle the prior for $K$ affects the Gibbs posterior distribution's properties.

\section{Gibbs posterior concentration rates}
\label{S:theory}

\subsection{Setup and objectives}
\label{SS:setup}

At a high level, our goal is to establish that the Gibbs posterior distribution, $\Pi_n$, for $\psi$ will concentrate its mass in vanishingly small neighborhoods of $\psi^\star$ as $n \to \infty$.  Then the radius of these vanishingly small neighborhoods would determine the {\em rate} at which the Gibbs posterior concentrates around $\psi^\star$.  More formally, we say that the Gibbs posterior {\em concentrates around a limit $\psi^\dagger$, with respect to a metric $d$, at rate $\eps_n \to 0$ if}
\[ \E \, \Pi_n(\{\psi: d(\psi^\dagger, \psi) > M_n \eps_n\}) \to 0, \quad n \to \infty, \]
where $M_n$ is a deterministic sequence that diverges to $\infty$ arbitrarily slowly.  The idea is that the Gibbs posterior is putting vanishingly small mass outside a neighborhood of $\psi^\dagger$ whose radius is $M_n \eps_n \approx \eps_n \to 0$.  

Throughout this section, we will focus on the case where $d$ is the usual $L_2(D)$ metric.  The particular rate that can be attained depends, in a specific way, on how the complexity index, $K$, of the model is being handled.  In particular, when $K$ is taken to be a fixed constant, so that $\tilde\pi(K)=1$ and there are no sums in \eqref{eq:gibbs}, Theorem~\ref{thm:K.fixed} states that the Gibbs posterior concentrates at a fast/optimal rate, $\eps_n = t_n^{-1/2}$, but around $\psi^\dagger=\psi_K^\perp$ instead of $\psi^\dagger=\psi^\star$.  Of course, if $\psi^\star \in \S_K$ for the chosen $K$, then $\psi_K^\perp = \psi^\star$ and we obtain the best possible result. This is not fully satisfactory, however, because generally one cannot expect $\psi^\star \in \S_K$ for a known $K$.  For this reason, we will consider two additional cases: one where we know the smoothness of $\psi^\star$ and can incorporate that information to select a complexity index sequence $K_n \to \infty$, and one where no such information is available and the posterior needs to learn the complexity index from the data itself.  For these two cases, the rate $\eps_n$ at which the Gibbs posterior concentrates around $\psi^\star$ will depend on the smoothness of $\psi^\star$ in a certain way (see below).  In the latter case, we say the rate is {\em adaptive} because the Gibbs posterior construction has no knowledge of the underlying smoothness of $\psi^\star$ and yet it can, in some sense, learn the appropriate complexity index from the data and concentrate at the same optimal rate as if that smoothness information were known.  See Theorems~\ref{thm:K.increasing}--\ref{thm:K.unknown} below for details.  

To end this first subsection, we will discuss the conditions we require on the sampling frequency, $\Delta=\Delta_n$, and how it depends on both $n$ and certain features of the basis functions chosen for the construction of the sieve $\{\S_K: K \geq 1\}$.  Recall that, for each $K \geq 1$, the collection $f_{K1},\ldots,f_{KK}$ are orthonormal $L_2(D)$ functions.  Here and throughout, we will also assume that they are each continuously differentiable on their respective supports.  For each $K$-specific collection, define the following two features:
\begin{align*}
F_1(K) & = \max_{k=1,\ldots,K} \Bigl\{ \sup_{x \in D} |f_{Kk}(x)| + \int_D |f_{Kk}'(x)| \, dx \Bigr\} \\
F_2(K) &= \max_{k=1,\ldots,K} \Bigl\{ \sup_{x \in D} |f^2_{Kk}(x)| + 2\int_D |f_{Kk}(x)f_{Kk}'(x)| \, dx \Bigr\}.
\end{align*}
Control on the magnitudes of these two features, as we discuss below, is effectively what 
\citet[p.~127--128]{figueroa2009nonparametric} states as his ``Standing Assumption~1.'' There are so many different choices of basis functions, but we have in mind the following two.  

\begin{enumerate}
\item {\em Piecewise polynomials}. 
The space of piecewise polynomials of degree at most $J-1$, based on a regular $L$-piece partition of window $D$ (i.e., $a=x_0<x_1<\cdots<x_L=b$), can be spanned by $J\times L$ orthonormal functions. If $\S_K$ denotes such a space, then its complexity index is $K=(J\times L)$, where the degree of polynomials is often taken as fixed, and 
\[ \S_K = \text{span}\{f_{K(j,l)}: j=0,\ldots,J-1, \ell=1,\ldots,L\}. \]
Specifically, basis functions are given as:
\begin{equation*}
    f_{K(j,l)}(x)=\bigl(\tfrac{2j+1}{x_{l}-x_{l-1}}\bigr)^{1/2} Q_j\bigl(\tfrac{2x-(x_{l-1}+x_{l})}{x_{l}-x_{l-1}}\bigr) \, 1_{(x_{l-1},x_{l})}(x),
\end{equation*}
where $Q_j$ is the Legendre polynomials defined on $[-1,1]$ of order $j$, satisfying 
\[ \int_{-1}^1 Q_j(x)Q_{j'}(x)dx=0, \quad \text{ for $j\neq j'$}. \]
\item {\em Trigonometric functions}. In this case, there is only one sequence of functions, rather than a collection of sequences indexed by $K$.  That is, the sieve is constructed as 
\begin{equation}
\label{eq:sieve}
\S_K = \text{span}\{f_1,\ldots,f_K\}, \quad K \geq 1. 
\end{equation}
What makes this formulation special is that the sieves are nested in the sense that $\S_K \subset \S_{K+1}$ for each $K$, which is important for the case considered in Section~\ref{SS:K.unknown}.  With $D=[a,b]$, the particular basis functions are defined as
\begin{align*}
f_1(x) & \equiv (b-a)^{-1/2}  \\
f_k(x) & = 
\begin{cases}
\bigl( \frac{2}{b-a} \bigr)^{1/2} \cos\bigl\{ \frac{k \pi (x-a)}{b-a} \}, & \text{if $k$ is even} \\
\bigl( \frac{2}{b-a} \bigr)^{1/2} \sin\bigl\{ \frac{(k-1) \pi (x-a)}{b-a} \}, & \text{if $k > 1$ is odd}.
\end{cases}
\end{align*}
\end{enumerate}
One key feature of both collections of functions is their approximation properties.  For example, \citet{MR1679028} states that, for Besov-smooth functions $\psi^\star$ with smoothness index $\alpha$, as defined below, the projection $\psi_K^\perp$ of $\psi^\star$ onto $\S_K$ satisfies 
\begin{equation}
\label{eq:approx}
\|\psi_K^\perp - \psi^\star\|_{L_2(D)} \lesssim K^{-\alpha}.
\end{equation}
Of relevance to the present discussion, note that, for piecewise polynomials, $F_1(K) \lesssim K^{1/2}$ and $F_2(K) \lesssim K$ and, for trigonometric functions, $F_1(K) \lesssim K$ and $F_2(K) \lesssim K$.  

Then the general requirement on the sampling frequency $\Delta$ can be described as follows.  Let $\kappa_n$ denote the largest index $K$ in consideration; in Sections~\ref{SS:K.fixed}--\ref{SS:K.unknown} this $\kappa$ will be a fixed constant, an increasing deterministic sequence, and an increasing deterministic upper bound on the support of the prior for $K$, respectively.  Using the fact that $t_n = n\Delta$, we write the requirement on the sampling frequency, $\Delta$, as  
\begin{equation}
\label{eq:Delta}
\max\{ F_1^2(\kappa_n) n \Delta^3, F_2(\kappa_n) \Delta\} = O(1), \quad n \to \infty. 
\end{equation}
In particular, when $K$ is a fixed constant, as in Theorems~\ref{thm:K.fixed}, condition \eqref{eq:Delta} boils down to $n\Delta^3=O(1)$ for any types of basis functions. In Theorem~\ref{thm:K.increasing}, where $K=K_n$ is an increasing, for both piecewise polynomials and trigonometric bases, condition \eqref{eq:Delta} is satisfied if $n \Delta^2 = O(1)$, assuming $\alpha>\frac12$. 
These constraints on $n$ and $\Delta$ are consistent with those in \citet{figueroa2009nonparametric,MR2787609} and \citet{adaptive.levy.2011}.  When $K$ is assigned a prior, as in Theorem~\ref{thm:K.unknown}, where we focus on the trigonometric basis functions, \eqref{eq:Delta} requires implies a slightly stronger condition, namely, $n \Delta^{5/3} = O(1)$.  Further comments on this latter point are given in the discussion following Theorem~\ref{thm:K.unknown}.



\subsection{Case 1: Fixed $K$}
\label{SS:K.fixed}

We start here with the simplest of cases, namely, where the model complexity index $K$ is taken to be fixed.  That is, there is no prior for $K$ so the Gibbs posterior for $\psi$ in \eqref{eq:gibbs} takes a simplified form, 
\[ \Pi_n(A) = \frac{\int_{\{\theta_K \in \Theta_K: \psi_{\theta_K} \in A\}} e^{-\omega t_n R_{n,K}(\theta_K)} \, \tilde\Pi_K(d\theta_K)}{\int_{\Theta_K} e^{-\omega t_n R_{n,K}(\theta_K)} \, \tilde\Pi_K(d\theta_K)}, \quad A \subseteq L_2(D). \]
This simplification makes computations relatively straightforward and, moreover, as we see in the theorem below, it allows for fast asymptotic concentration of the Gibbs posterior.  However, the same simplification is a restriction on the model's flexibility, so there is generally a bias that cannot be overcome, even asymptotically.  Indeed, the result below states that the Gibbs posterior concentrates quickly around $\theta_K^\perp$, not around $\psi^\star$.  Of course, if the model is ``correct'' in the sense that $\psi^\star \in \S_K$, then there is no bias and we achieve the best possible concentration rate result.  

\begin{cond}
\label{cond:K.fixed}
\mbox{}
\begin{enumerate}
\item $\Delta$ satisfies \eqref{eq:Delta} for $\kappa_n \equiv K$.
\item The $K$-specific prior distribution, $\tilde\Pi_K$, has a density function that is continuous and bounded away from 0 in a neighborhood of $\theta_K^\perp$. 
\end{enumerate}
\end{cond}

This is a standard condition in the literature on Bayesian asymptotics and would hold for virtually any choice of prior for the $K$-vector $\theta_K$.  

\begin{thm}
\label{thm:K.fixed}
Under Condition~\ref{cond:K.fixed}, the Gibbs posterior distribution, $\Pi_n$, satisfies 
\[ \E \, \Pi_n(\{\psi: \|\psi - \psi_K^\perp\|_{L_2(D)} > M_n \eps_n\}) \to 0, \quad n \to \infty, \]
where $\eps_n = t_n^{-1/2}$ and $M_n \to \infty$ arbitrarily slowly.  In particular, if $\psi^\star \in \S_K$, then the Gibbs posterior concentrates around $\psi^\star$ at rate $\eps_n = t_n^{-1/2}$.  
\end{thm}

Since $\Delta$ is vanishing, the rate $\eps_n = t_n^{-1/2}$ is slower than what is often achieved in regular, finite-dimensional estimation problems.  This slower rate is a consequence of the problem, not our proof technique.  Indeed, the accumulation of information in this problem is not linear in $n$, since $n \to n+1$ provides an additional observation but also extends the time horizon over which observations are taken.  So while the $(n+1)^\text{st}$ observation provides extra information, a small amount of information is lost as a result of extending the time horizon.

\subsection{Case 2: Increasing $K=K_n$}
\label{SS:K.increase}

To achieve the fast concentration rate in Theorem~\ref{thm:K.fixed} to the true L\'evy density, $\psi^\star$, requires the strong assumption that $\psi^\star \in \S_K$.  In words, it requires $\psi^\star$ to be a very nice function.  We may not always be willing to make such a strong assumption, so more flexibility needs to be built in to the model.  The additional flexibility we consider here is to allow the model complexity index $K$ be increasing in $n$.  That is, the fixed $K$ appearing in the definition of the Gibbs posterior, $\Pi_n$, for $\psi$ in  Section~\ref{SS:K.fixed} is replaced by an increasing sequence $K_n$.  To determine how fast $K=K_n$ should be growing in order to achieve the desired flexibility, we need to be more specific about our assumption on $\psi^\star$.  

As is customary in the nonparametric estimation literature, here we will assume that $\psi^\star$, or at least its restriction, $\psi^\star 1_D$, to $D$, is smooth in the sense that it belongs to a certain Besov class of functions \citep{MR1679028}.  For a function $g: D \to \mathbb{R}$, define the first-order difference operator 
\[ \delta_h^1(g, x) = g(x + h) - g(x), \quad h > 0, \]
and then recursively define the higher-order differences 
\[ \delta_h^r(g, x) = \delta_h^1\{ \delta_h^{r-1}(g,\cdot), x\}, \quad \text{$r$ a positive integer}. \]
Then we say that the function $g$ belongs to the Besov space $\mathcal{B}_{\infty}^{\alpha} = \mathcal{B}_{\infty}^{\alpha}(L_p(D))$ with $p\geq2$ and $\alpha>0$, if  $\|g\|_{\mathcal{B}_{\infty}^{\alpha}}<\infty$, where 
\[ \|g\|_{\mathcal{B}_{\infty}^{\alpha}} := \sup_{\gamma > 0} \gamma^{-\alpha} \sup_{h \in (\gamma, 0)} \| \delta_h^r(g, \cdot)\|_{L_p(D)}, \quad r = [\alpha] + 1, \]
with $[\alpha]$ denoting the integer part of $\alpha$.  The index $\alpha$ determines the smoothness of the functions in $\mathcal{B}_{\infty}^{\alpha}$ and will be important for determining the (optimal) Gibbs posterior concentration rate; see below. What follows holds for all fixed $p \geq 2$. 

This boils down to a nonparametric function estimation problem so we can anticipate what would be the optimal rate we hope to achieve and how it depends on the smoothness of $\psi^\star$.  \citet{figueroa2009nonparametric} showed that the minimax optimal rate for the L\'evy density estimation problem in Besov spaces is 
\[ \inf_{\hat\psi} \sup_{\psi^\star \in \mathcal{B}_\alpha^\infty} \E\|\hat\psi - \psi^\star\|_{L_2(D)}^2 \gtrsim t_n^{-2\alpha/(2\alpha + 1)}. \]
That is, for the ``most difficult'' of $\psi^\star$ to estimate, the best possible estimator would have a rate like that in the right-hand side of the above display.  So, if we can show that our Gibbs posterior concentrates around $\psi^\star$ at this rate, then we will know that it is minimax optimal.  The particular choice of $K=K_n$ is made by setting the approximation error bound on the right-hand side of \eqref{eq:approx} equal to the minimax optimal rate, $t_n^{-\alpha/(2\alpha+1)}$, and solving for $K$.  It is easy to see that the solution, $K_n$, is roughly $t_n^{1/(2\alpha+1)}$, which is the choice taken Theorem~\ref{thm:K.increasing} below.   

\begin{cond}
\label{cond:K.increasing}
\mbox{}
\begin{enumerate}
\item $\psi^\star 1_D$ is bounded and $\psi^\star$ belongs to the Besov class $\mathcal{B}_{\infty}^{\alpha}$, with $\alpha$ known.  
\item For the same $\alpha$, let $K_n = \ceiling{t_n^{1/(2\alpha+1)}}$, where $\ceiling{z}$ denotes the smallest integer greater than or equal to $z$.  Then $\Delta$ satisfies \eqref{eq:Delta} with $\kappa_n = K_n$. 
\item For the prior, $\tilde\Pi_K$, let the components of $\theta_{Kk}$, $k=1,\ldots,K$, be independent with a positive density function at $\theta_{Kk}^\perp$. 
\end{enumerate}
\end{cond}


\begin{thm}
\label{thm:K.increasing}
Under Condition~\ref{cond:K.increasing}, the Gibbs posterior distribution $\Pi_n$ satisfies 
\[ \E\,\Pi_n(\{\psi: \|\psi - \psi^\star\|_{L_2(D)} > M_n \eps_n\}) \to 0, \quad n \to \infty, \]
where $\eps_n = t_n^{-\alpha/(2\alpha+1)}$ and $M_n \to \infty$ arbitrarily slowly. 
\end{thm}

As discussed above, the rate achieved by the Gibbs posterior is minimax optimal.  It is also worth pointing out that this is the same optimal rate that could be achieved if the entire sample path $X$ were observable, i.e., if we had continuous sampling, which is quite remarkable.  Furthermore, there is no extra logarthmic term in the rate achieved by the Gibbs posterior, as is often the case in nonparametric problems.  This is due to our special requirement on the prior $\tilde\Pi_{K_n}$.  If we relax the independence assumption, then a version of Theorem~\ref{thm:K.increasing} still holds, but with a slightly slower rate $\eps_n' = (\log t_n)^{1/2} t_n^{-\alpha/(2\alpha+1)}$.  

Finally, note that Theorem~\ref{thm:K.fixed} is effectively a special case of Theorem~\ref{thm:K.increasing}.  That is, if $\alpha \to \infty$, so that $\psi^\star$ is a ``nice'' function, in the sense of being very smooth (e.g., in $\S_K$ for some $K$), then $K_n$ is effectively a constant and the rate achieved is roughly $t_n^{-1/2}$.

\subsection{Case 3: Prior distribution for $K$}
\label{SS:K.unknown}

The only disadvantage to the result in Theorem~\ref{thm:K.increasing} is that it requires knowledge of the smoothness $\alpha$ of $\psi^\star$.  Generally, this $\alpha$ would be unknown and, therefore, it would be impossible to set $K_n = \ceiling{t_n^{1/(2\alpha+1)}}$ in the Gibbs construction.  Alternatively, as we described above, instead of taking $K$ fixed and associating $\psi$ with a $K$-vector $\theta_K$, we can express $\psi$ as the pair $(K,\theta_K)$ and introduce a prior for the pair.  Recall that the marginal prior for $K$ is denoted as $\tilde\pi$.  Provided that $\tilde\pi$ assigns its mass in a suitable way, we can show that (roughly) the same optimal Gibbs posterior concentration rate as in Theorem~\ref{thm:K.increasing} is achieved, {\em adaptively}, without prior knowledge of $\alpha$.  

Compared to the previous sections, here we are simultaneously considering sieves indexed by different $K$ and our theoretical analysis requires that there be some meaningful connection between the spaces $\S_K$ for different $K$.  A natural connection is for these spaces to be nested in the sense that $\S_K \subseteq \S_{K+1}$, $K \geq 1$.  To achieve this, we take a concrete choice of sieves, namely that of the form \eqref{eq:sieve} with the functions $(f_k)$ being the trigonometric functions as described in Section~\ref{SS:setup} above.  

The prior for $(K,\theta_K)$ is specified hierarchically.  The marginal prior for $K$ is 
\begin{equation}
\label{eq:prior.for.K}
\tilde\pi(K) \propto e^{-\beta K \log K}, \quad K=1,\ldots, \ceiling{t_n}, 
\end{equation}
where $\beta > 0$ is a constant to be specified by the data analyst; see below.  Note that the marginal prior for $K$ is truncated so that $K > t_n$ has prior probability 0.  This is reasonable because (a)~Occam's razor suggests that we should attempt to limit the complexity of the model being fit, and (b)~we know that the ``ideal'' $K$ is $\ceiling{t_n^{1/(2\alpha+1)}}$, which is never more than $t_n$.  For the conditional prior of $\theta_K$, given $K$, we make the same assumption as in Theorem~\ref{thm:K.increasing}.  

\begin{cond}
\label{cond:K.unknown}
\mbox{}
\begin{enumerate}
\item $\psi^\star 1_D$ is bounded and belongs to the Besov class $\mathcal{B}_{\infty}^{\alpha}$ with $\alpha$ unknown. 
\item The (nested) sieves $(\S_K)$ are of the form \eqref{eq:sieve}, where the $f_K$'s are the trigonometric basis functions.  
\item The condition \eqref{eq:Delta} holds with $\kappa_n = \ceiling{t_n}$ which, for the trigonometric basis functions, requires that $n \Delta^{5/3} = O(1)$. 
\item For each $K \geq 1$, the conditional prior $\tilde\Pi_K$ assigns independent prior distributions to each basis coefficient $\theta_k$, for $k=1,\ldots,K$, each having positive density function in a neighborhood of $\theta_k^\perp$. 
\item The marginal prior, $\tilde\pi$, for $K$ is of the form \eqref{eq:prior.for.K}, where the coefficient $\beta > 0$ satisfies
\begin{equation}
\label{eq:beta}
\beta > \omega C^2 + o(1),
\end{equation}
with $C > 0$ given by
\[ C^2 > 2 \sup_{x \in D} \psi^\star(x) + o(1), \quad n \to \infty. \]
The above condition on $\beta$ can equivalently be written as a requirement that the learning rate $\omega$ be sufficiently small in the sense that 
\begin{equation}
\label{eq:omega}
\omega < C^{-2} \{ \beta + o(1)\}. 
\end{equation}
\end{enumerate}
\end{cond}

\begin{thm}
\label{thm:K.unknown}
Under Condition~\ref{cond:K.unknown}, the Gibbs posterior distribution satisfies 
\[ \E\,\Pi_n(\{\psi: \|\psi - \psi^\star\|_{L_2(D)} > M_n \eps_n\}) \to 0, \quad n \to \infty, \]
where $\eps_n = (\log t_n)^{1/2} t_n^{-\alpha/(2\alpha + 1)}$ and $M_n \to \infty$ arbitrarily slowly.  
\end{thm}

Recall that the difference between Theorems~\ref{thm:K.increasing} and \ref{thm:K.unknown} is that the latter achieves (nearly) the minimax optimal rate without the prior knowledge of $\alpha$ that the former assumed.  So the extra logarithmic term in Theorem~\ref{thm:K.unknown}'s rate can be interpreted as the (negligible) cost of adaptation to the unknown smoothness of $\psi^\star$.  

Compared to the case where the smoothness is known, here we have a slightly tighter constraint on $\Delta$, coming from the stronger condition $n\Delta^{5/3} = O(1)$.  \citet{figueroa2009nonparametric} did not get an adaptive convergence rate result under the weaker $n\Delta^2 = O(1)$ condition.  \citet{adaptive.levy.2011}, on the other hand, do achieve adaptive convergence rates for point projection estimators under that condition.  It is unclear if the difference is due to our proof techniques or due to our choice of basis.  Regardless, our conclusion in Theorem~\ref{thm:K.unknown} above is arguably still stronger in the sense that we obtain a rate for a full probability distribution on the space of L\'evy densities, as opposed to just a point estimator.  Moreover, by averaging over different complexity indices according to a marginal Gibbs posterior, we are better accounting for the uncertainty in $K$ compared to previous approaches that simply plug in a fixed choice.  

We can also assess the Gibbs posterior's ability to learn the complexity of $\psi^\star$ by looking at the asymptotic properties of the marginal Gibbs posterior for $K$.  In particular, a relevant question is if the Gibbs posterior tends to ``overfit'' in the sense of supporting values of $K$ that are much larger than the complexity index $K_n$ used by an oracle, like in the context of Theorem~\ref{thm:K.increasing}, who knows the smoothness $\alpha$ of $\psi^\star$.  The following result shows that, under conditions comparable to those in Theorem~\ref{thm:K.unknown}, the marginal Gibbs posterior assigns vanishing probability to $K$ values larger than a constant multiple of $K_n$.

\begin{thm}
\label{thm:no.overfit}
Suppose Condition~\ref{cond:K.unknown} holds, but with \eqref{eq:beta} replaced by 
\begin{equation}
\label{eq:beta.new}
\beta > \tau(\tau - 1)^{-1} C^2 \omega, \quad \text{for some $\tau > 1$}. 
\end{equation}
Then $\E\,\Pi_n(\{K: K > \tau K_n\}) = o(1)$ as $n \to \infty$. 
\end{thm}

Note that, to achieve the strongest conclusion from Theorem~\ref{thm:no.overfit}, the condition \eqref{eq:beta.new} needs to be a bit stronger than the condition \eqref{eq:beta} used in the proof of Theorem~\ref{thm:K.unknown}.  That is, the conclusion is strongest for $\tau$ just slightly larger than 1, in which case, the right-hand side of \eqref{eq:beta.new} would be larger than the right-hand side of \eqref{eq:beta}.  

It is perhaps not surprising, that the conditions on $\beta$ stated in Theorems~\ref{thm:K.unknown}--\ref{thm:no.overfit} depend on a certain feature of the unknown $\psi^\star$, namely, $\sup_{x \in D} \psi^\star(x)$.  Since the learning rate $\omega$ can be taken to be arbitrarily small without affecting the concentration rate, the conditions \eqref{eq:beta} or \eqref{eq:beta.new} can always be arranged so, in theory, there is no constraint.  In practice, however, the user must set specific values for $(\beta,\omega)$, and likely this choice would need to be data-driven.  General methods for selecting suitable values of $(\omega,\beta)$ are beyond the scope of the present paper, but one simple idea is suggested in Section~\ref{S:numerical}.



\section{Practical implementation}
\label{S:numerical}

\subsection{Posterior computation}

In addition to being free of a model/likelihood and achieving the optimal posterior concentration rates, important feature of our proposed framework for inference on a general L\'evy density is that computation is not only doable, it is also relatively straightforward and efficient. Below we present these details for the known- and unknown-$K$ cases.  


In the case where $K$ is known, as discussed in Section~\ref{SS:K.fixed} and Section~\ref{SS:K.increase}, the Gibbs posterior can be considered as a distribution for the $K$-dimensional basis coefficient $\theta_K$. The term  $e^{-\omega t_n R_{n,K}(\theta_K)}$ in \eqref{eq:gibbs} is proportional to the density of a multivariate normal distribution for $\theta_K$, with independent components.  That is, 
\[ e^{-\omega t_n R_{n,K}(\theta_K)} = \prod_{k=1}^K \nm\bigl(\theta_{Kk} \mid \hat\theta_{Kk}, (2\omega t_n)^{-1} \bigr). \]
Therefore, if we choose the prior $\tilde\Pi_K$ for $\theta_K$ as 
\begin{equation}
\label{eq:iid.prior}
\theta_{Kk} \iid \nm(0, \sigma_0^2), \quad k=1,\ldots,K, 
\end{equation}
then the Gibbs posterior for $\theta_K$ satisfies
\[ \theta_{Kk} \ind \nm\big(\hat\theta_{Kk}(1+(2\omega t_n)^{-1}\sigma_0^{-2})^{-1},(2\omega t_n+\sigma_0^{-2})^{-1}\big), \quad k=1\cdots K. \]
Then by drawing posterior samples $\theta_K$ from such a multivariable normal distribution, posterior inference on the L\'evy density readily follows.  Our choice to use a normal prior---which is conjugate to the Gibbs ``likelihood''---means we do not need {\em Markov chain} Monte Carlo, just the ordinary independent Monte Carlo is enough.  Of course, other choices of prior can be considered, and the formulas above can be adjusted accordingly, but then sampling from the Gibbs posterior might be more expensive.  

In the case where $K$ is unknown, as discussed in Section~\ref{SS:K.unknown}, the Gibbs posterior is for the pair $(K,\theta_K)$. Therefore, to sample from the Gibbs posterior, we first sample $K$ from its marginal distribution, and then generate $\theta_K$ according to the conditional distribution on $K$. Specifically, if we use the same iid prior for the $\theta_k$'s as is \eqref{eq:iid.prior}, then the marginal Gibbs posterior distribution of $K$ has a probability mass function 
\[ \pi_n(K) \propto \exp\Big\{\sum_{k=1}^K\frac{\omega t_n\hat{\theta}^2_k}{1+(2\omega t_n)^{-1}\sigma_{0}^{-2}}-
\frac{K}{2} \log(2\omega t_n\sigma_{0}^2+1)
-\beta K \log K\Big\}. \]
The conditional Gibbs posterior distribution of $\theta_K$, given $K$, is exactly the same as that in the known $K$ case discussed above.

\subsection{Illustration}

Here we present only a brief numerical illustration to demonstrate that the proposed Gibbs posterior can be computed efficiently and provides an accurate estimation of the underlying L\'evy density as the results in Section~\ref{S:theory} suggest, even in large-but-finite sample cases. We consider only the unknown-$K$ case since that is the most practical. 

Following \citet{figueroa2009nonparametric}, we consider $X$ a so-called variance--gamma process, an important special case of the general L\'evy processes, as proposed by \citet{madan1998variance}. 
A variance--gamma process can be interpreted as Brownian motion with drift with a random time generated from a gamma process. That is, 
$X(t)=B(U(t;\nu);\mu,\sigma)$
where $B(t;\mu,\sigma)$ represents a
Brownian motion with a linear drift of slope $\mu$ and a volatility coefficient $\sigma$, and $U(t;\nu)$ represents a gamma process whose evaluation at time $t$ follows a gamma distribution with a shape parameter $t \nu^{-1}$ and a scale parameter $\nu$. Based on such an interpretation, discrete observations $X(t_i), i=0,\cdots, n$ with $t_i=i\Delta$ can be obtained by alternatively generating iid increments $Y_i=X(t_i)-X(t_{i-1}),i=0,\cdots n$, given as
\begin{align*}
U_i &\iid \mathsf{Gamma}(\Delta \nu^{-1},\nu),  \\
(Y_i \mid U_i) & \ind \nm(\mu U_i,\sigma^2 U_i), \qquad i=1,\ldots,n, 
\end{align*}
where $\mathsf{Gamma}(a,b)$ represents a gamma distribution with density function proportional to $x \mapsto x^{a-1} e^{-x/b}$.  
The L\'evy density corresponding to a variance--gamma process parameterized by $(\mu, \sigma, \nu)$ has the following expression
\begin{align*}
    \psi^\star(x)=\begin{cases}
      \nu^{-1}|x|^{-1}e^{x/\eta^{+}},&x>0\\
      \nu^{-1}|x|^{-1}e^{x/\eta^{-}},&x<0
    \end{cases}
\end{align*}
where $\eta^{\pm}=(\mu^2\nu^2/4+\sigma^2\nu/2)^{1/2} \pm \mu\nu/2$. In the following simulation, 
the triplet $(\mu, \sigma, \nu)$ is taken as $(0, 3.7\times10^{-1.5}, 2\times10^{-3})$, so we consider the same variance--gamma process set up as in \citet{figueroa2009nonparametric}. 

\begin{figure}
\begin{center}
\scalebox{0.7}{\includegraphics[scale=0.6]{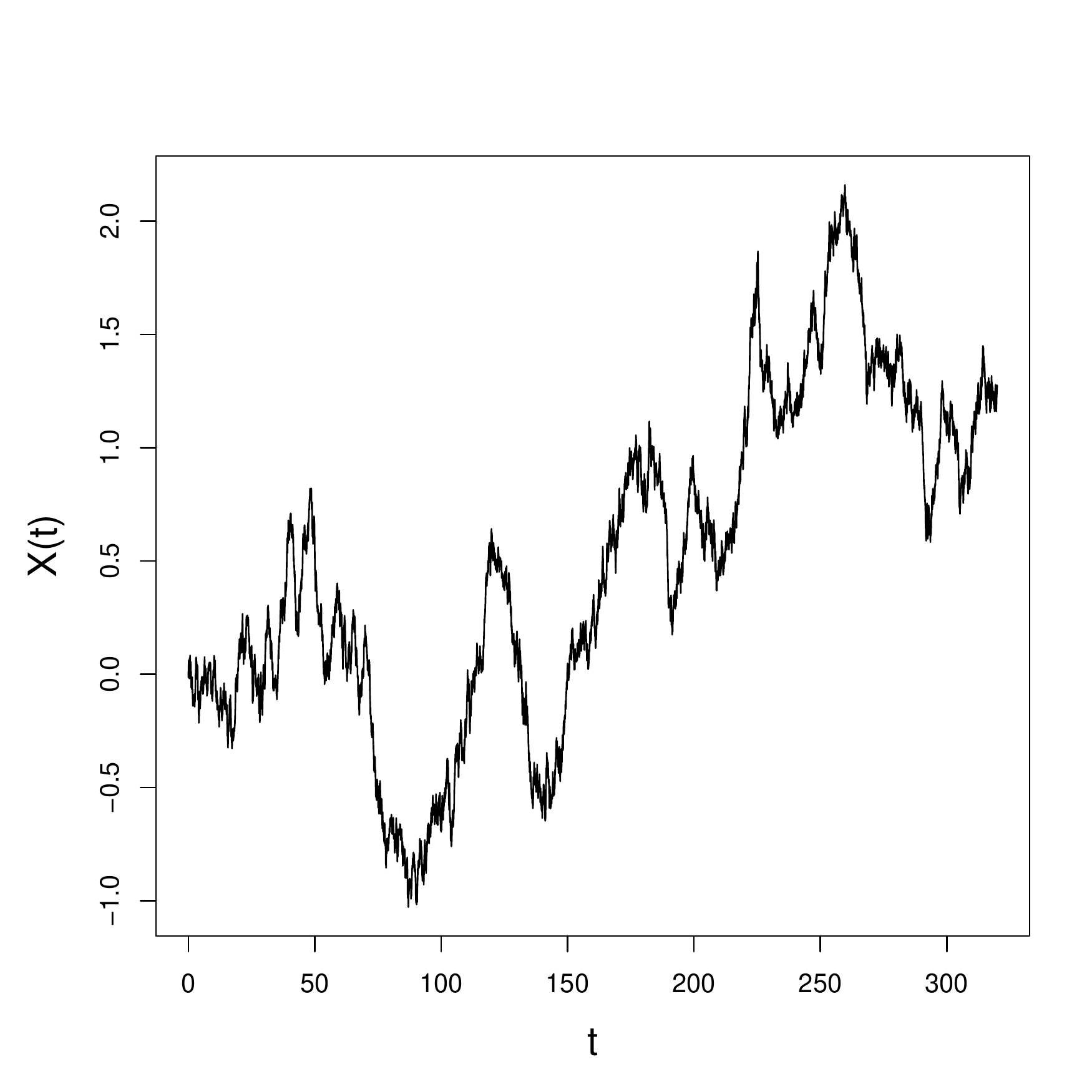}}
\end{center}
\caption{Plot of the discretely observed sample path $t \mapsto X(t)$, on the time window $[0,320]$, from the variance--gamma process described in the text.}
\label{fig:data}
\end{figure}

Here we consider three sampling settings where the number of observations, $n$, satisfies $n= \ceiling{0.05\times \Delta^{-5/3}}$, where $\Delta = 10^{-3} \times 2^{-3j}$, for $j=1,2,3$.   
As the index $j$ increases, we not only have more observations, but the spacings between those observations are shrinking, as one expects in a high-frequency sampling setting.  
Figure~\ref{fig:data} shows a plot of the (discretely-sampled version of the) sample path $t \mapsto X(t)$, which corresponds to the $j=3$ setting, where the last observation time is $t_n = n\Delta = 320$.  
Recall that we cannot separate the continuous and jump parts of the process without a high sampling frequency.  Moreover, the theory in the previous section implies that accurate estimation also requires a large number of observations.  So our choice of large $n$ and small $\Delta$ are consistent with what both intuition and our theory suggest is needed.  

To construct the Gibbs posterior, the prior distribution for $K$ is specified in \eqref{eq:prior.for.K} with $\beta=0.5$, and the conditional prior distribution for $\theta_K$, given $K$, is specified in \eqref{eq:iid.prior} with $\sigma_0=10^3$. Throughout, we use a fixed learning rate $\omega=10^{-5}$ and seek to estimate the L\'evy density on the window $D=[0.006, 0.014]$.  It is convenient to actually do the numerical calculations on a slightly wider interval $D'$, and then restrict those results to the target interval $D$; here we used $D'=[0.005, 0.015]$.  In a simulation study like this, where $\sup_{x \in D} \psi^\star(x)$ is known, we can set $\omega$ and $\beta$ to satisfy \eqref{eq:beta}. When the true L\'evy density is unknown, a naive strategy for choosing $(\omega,\beta)$ is 
to use \eqref{eq:beta} but with an estimator $\sup_{x \in D} \hat\psi(x)$ in place of that unknown feature of $\psi^\star$.  Whether this is an effective strategy remains an open question. 

\begin{figure}[t]
\begin{center}
\subfigure[$j=1$]{\scalebox{0.28}{\includegraphics{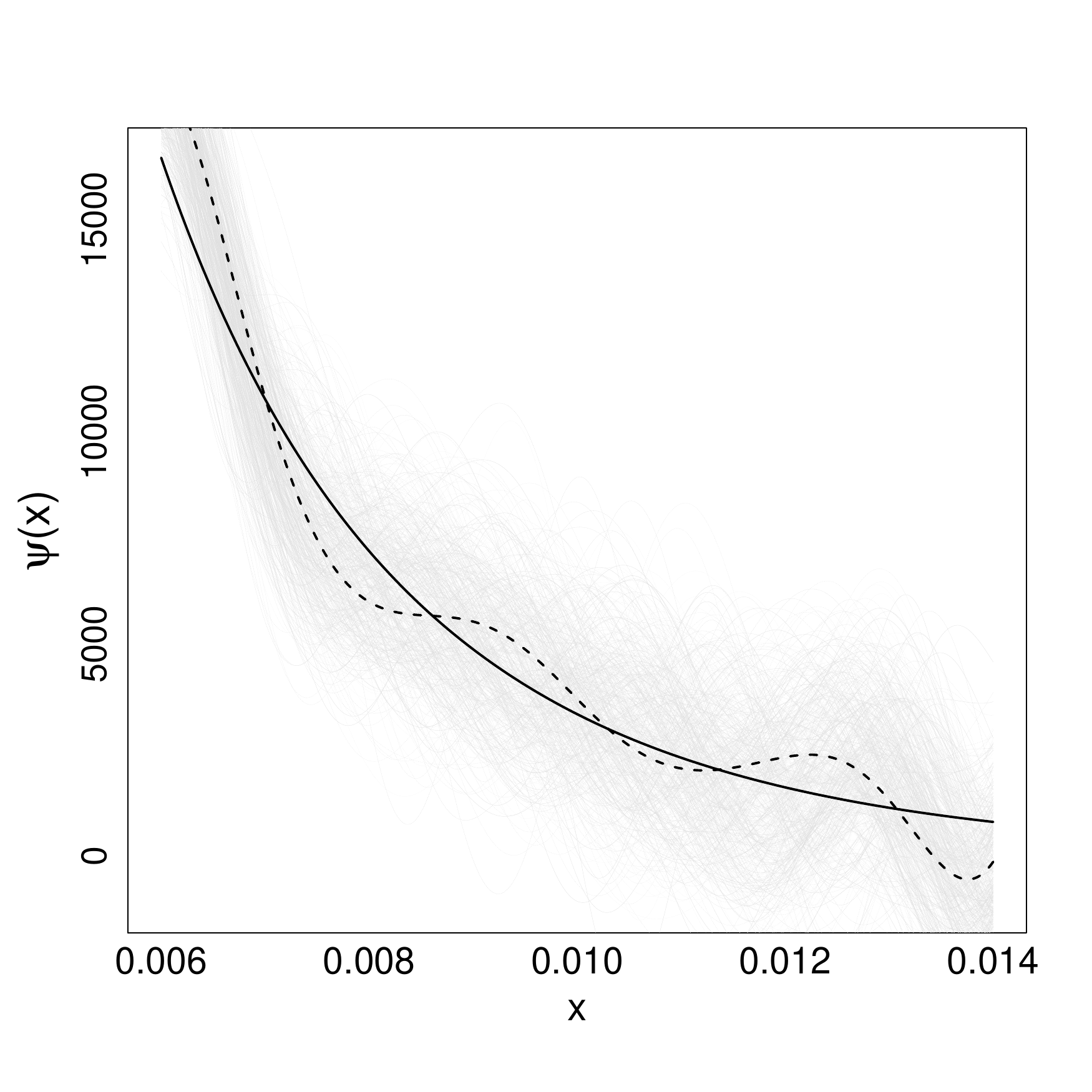}}}
\subfigure[$j=2$]{\scalebox{0.28}{\includegraphics{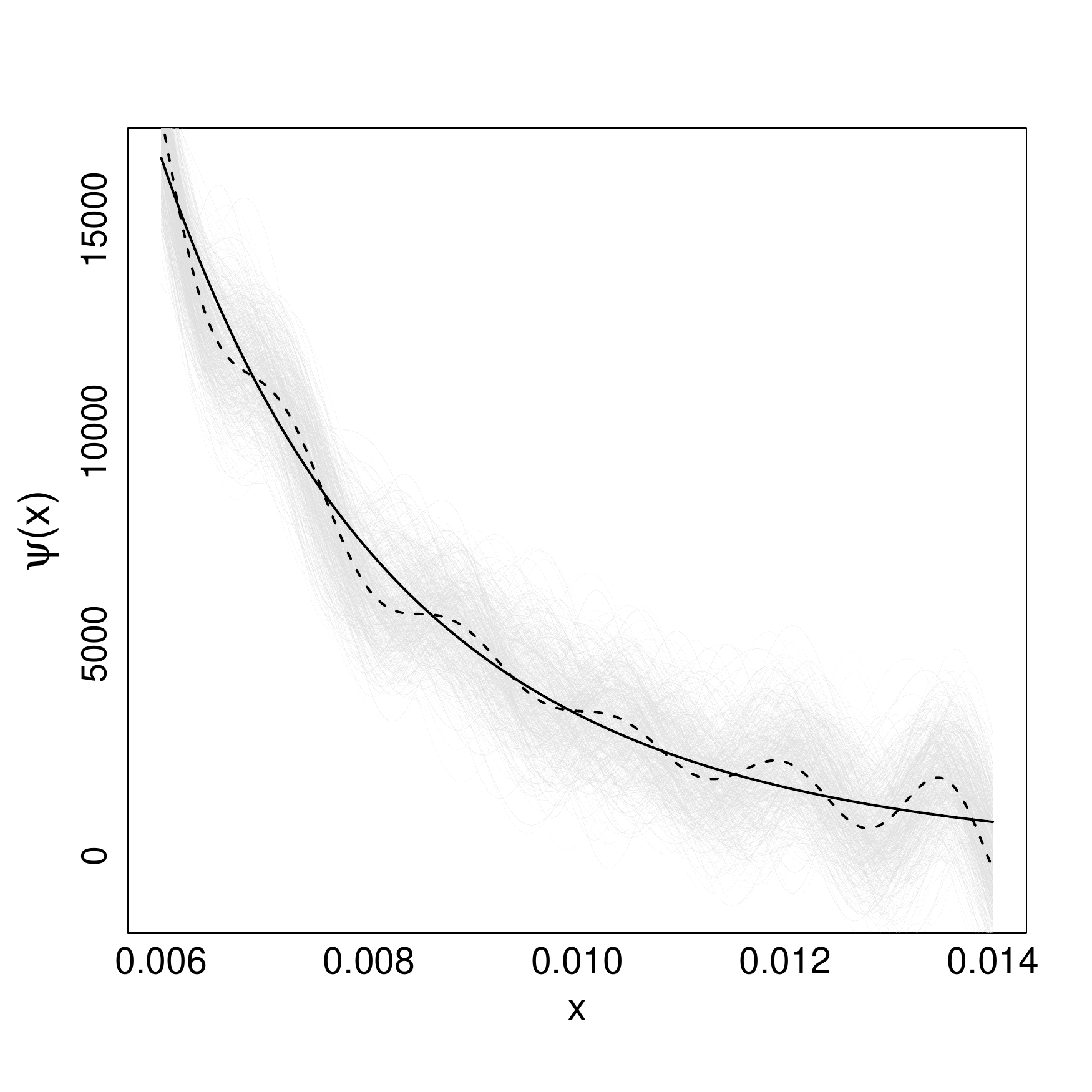}}}
\subfigure[$j=3$]{\scalebox{0.28}{\includegraphics{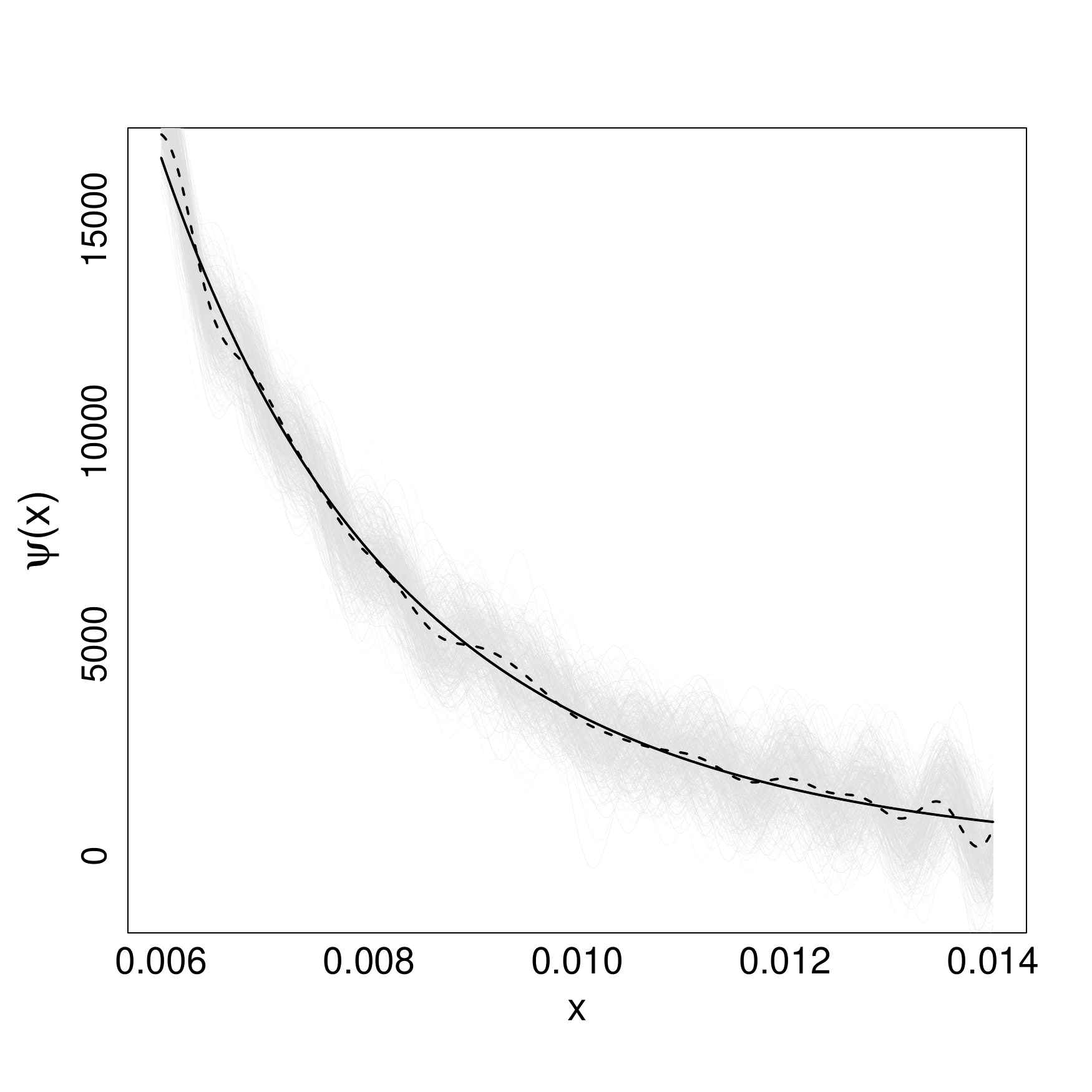}}}
\end{center}
\caption{Posterior samples of $\psi$ (gray), posterior mean function (dotted), and the true L\'evy density $\psi^\star$ (solid black), across the three sampling schemes indexed by $j$.}
\label{fig:post_psi}
\end{figure}

\begin{figure}[t]
\begin{center}
\subfigure[$j=1$]{\scalebox{0.28}{\includegraphics{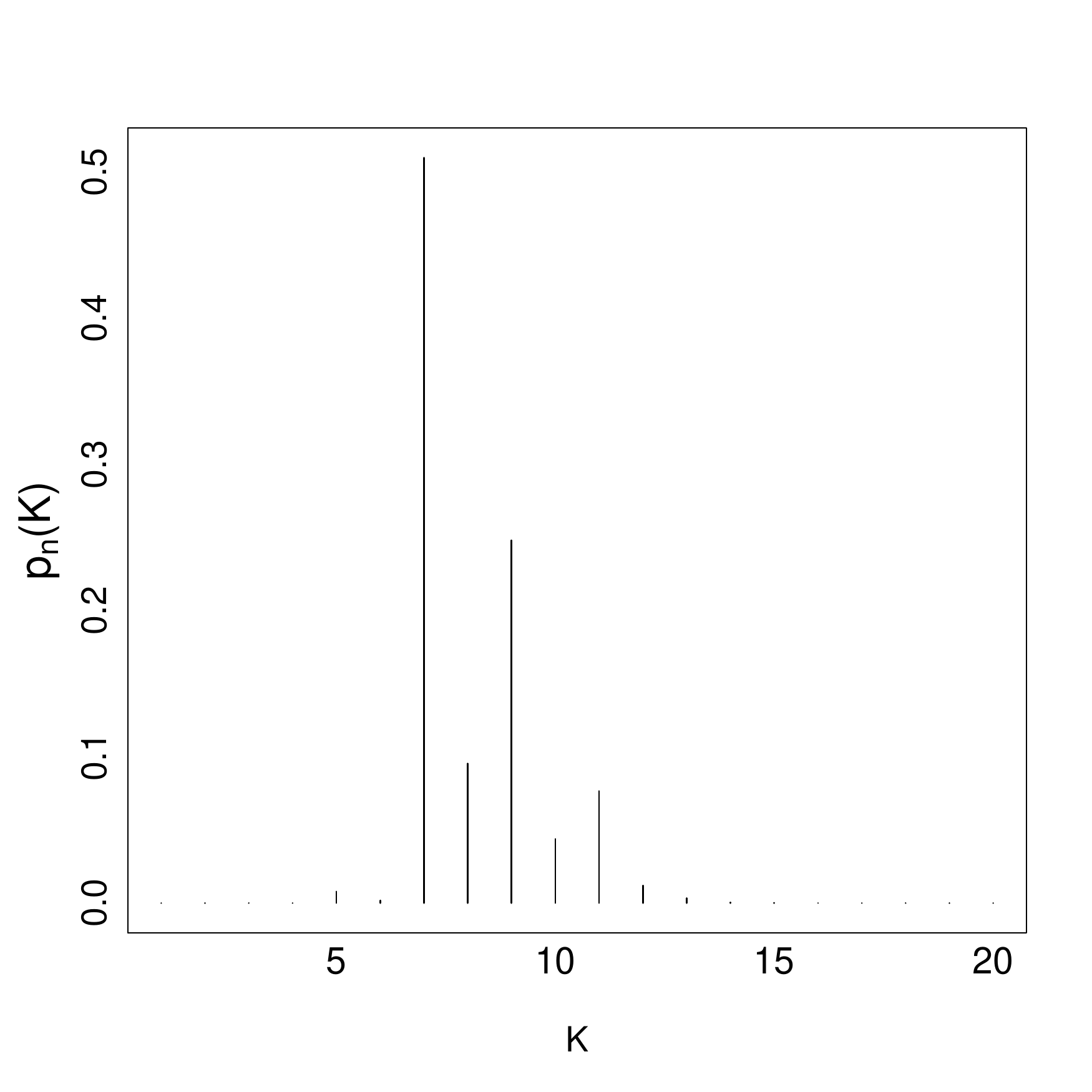}}}
\subfigure[$j=2$]{\scalebox{0.28}{\includegraphics{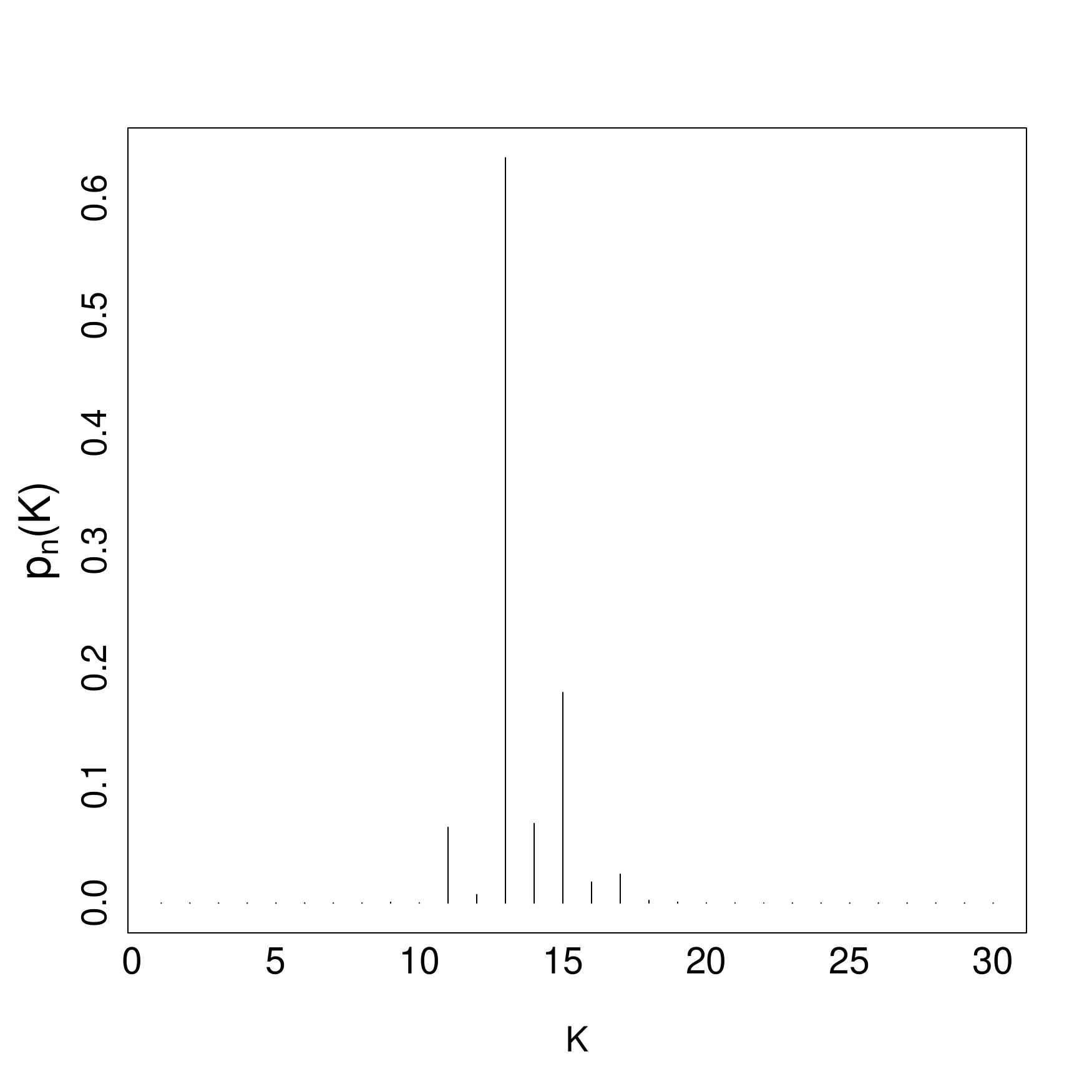}}}
\subfigure[$j=3$]{\scalebox{0.28}{\includegraphics{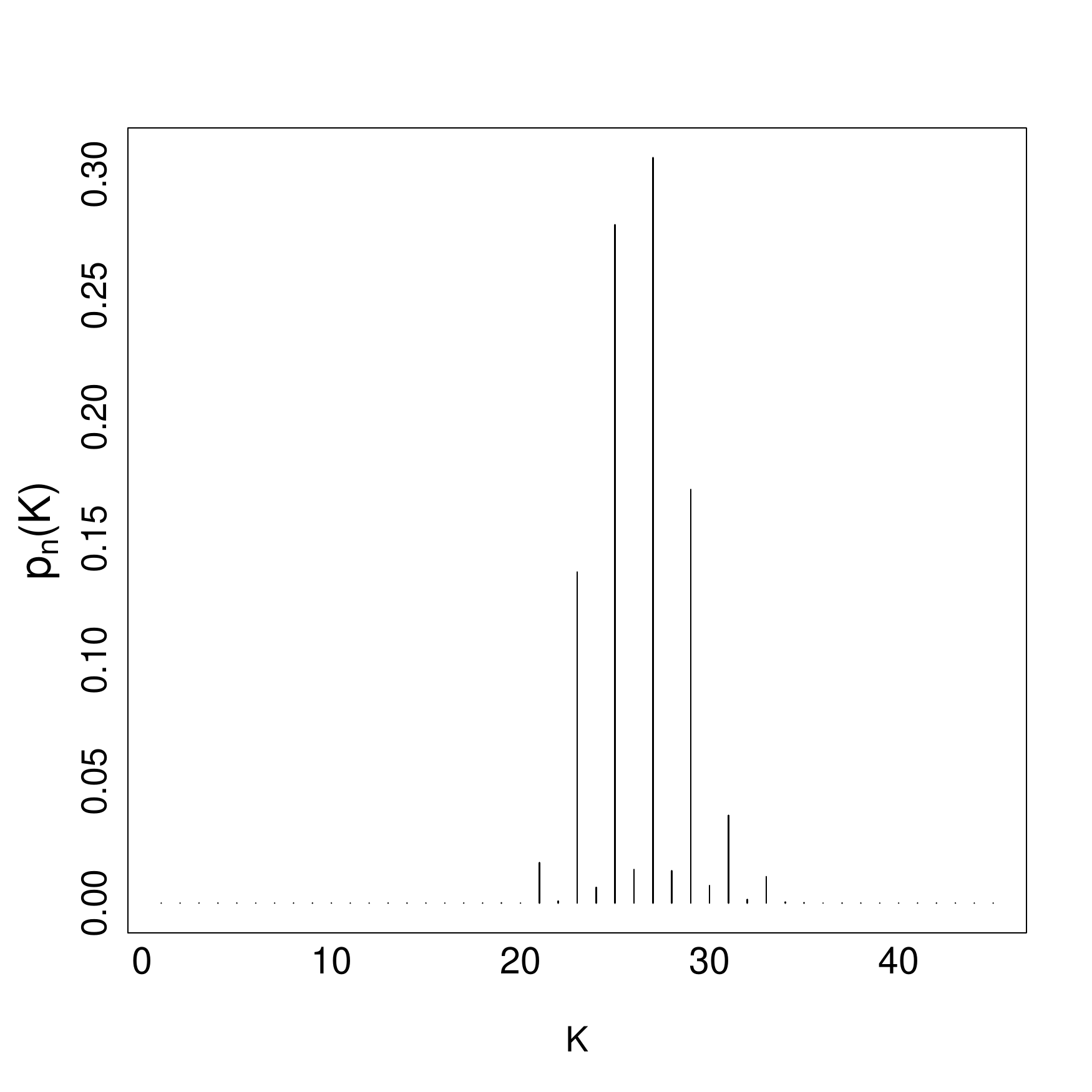}}}
\end{center}
\caption{Posterior marginal distribution of $K$, across the three sampling schemes.}
\label{fig:post_K}
\end{figure}

Applying our Gibbs posterior framework to the data presented in Figure~\ref{fig:data} produces the following results. 
Figure~\ref{fig:post_psi} shows 500 (of the total $1000$) Gibbs posterior samples of $\psi$ and the corresponding Gibbs posterior mean function, with the true L\'evy density overlaid.  Clearly, as $j$ increases---i.e., increasing $n$ and decreasing $\Delta$---the posterior samples and the posterior mean function becomes more concentrated around the true L\'evy density, as Theorem~\ref{thm:K.unknown} predicts, and that the posterior mean curve is smoother than any individual sample thanks to the mixing over the posterior distribution of $K$. Indeed, Figure~\ref{fig:post_K} shows the posterior marginal distributions of $K$, which as $j$ increases, concentrates its mass on larger values of the complexity index, which is consistent with our expectations based on Theorem~\ref{thm:no.overfit}; some further comments on this are given in Section~\ref{S:discuss}.  



\section{Discussion}
\label{S:discuss}

In this paper, we considered the challenging problem of estimation and inference on the L\'evy density under discrete sampling from a general L\'evy process.  Motivated by the inability to write down an explicit likelihood in the most general of cases, we leveraged an ``(approximate) expected loss minimizer'' interpretation of the L\'evy density to construct a so-called Gibbs posterior distribution for $\psi$.  Besides being model- or likelihood-free, this approach also allows for the incorporation of available prior information, e.g., on how smooth one might expect the true L\'evy density to be.  We investigated the asymptotic concentration properties of our proposed Gibbs posterior distribution and established, in the most interesting, challenging, and practically relevant cases, the posterior concentrates around the true $\psi^\star$ adaptively and at (nearly) the minimax optimal rate.  To our knowledge, this is the first Bayesian-like concentration rate result for this general L\'evy density estimation problem.  

Beyond the theory, it turns out that there are some practical advantages to this Gibbs posterior approach.  One is that, since it directly focuses on the L\'evy density, there are no nuisance parameters to assign prior distributions for and to marginalize out in the posterior computations.  Another advantage is that these posterior computations can be done in a straightforward and efficient manner, even in the case of a hierarchical model for $(K,\psi_K)$ as discussed in Section~\ref{SS:K.unknown}.  And the brief numerical illustration in Section~\ref{S:numerical} highlights this simplicity and demonstrates that the good performance suggested by the asymptotic concentration rate results carries over to finite-sample situations as well.  

While the asymptotic concentration results hold for any fixed choice of $\omega$ (with an appropriate corresponding choice of $\beta$ in the context of Theorem~\ref{thm:K.unknown}), the fact that $\omega$ does affect the finite-sample performance means that data-driven choices are needed.  As we briefly mentioned in Section~\ref{S:gibbs}, there is an active literature on data-dependent learning rate selection; see \citet{gpc.compare} for a comparison.  However, it remains unclear how these various methods might perform in an especially challenging application such as this, especially one with an infinite-dimensional quantity of interest.  

Another interesting open question concerns the asymptotic concentration properties of the marginal Gibbs posterior for $K$ in the context of Theorem~\ref{thm:K.unknown}.  We showed in Theorem~\ref{thm:no.overfit} that this marginal Gibbs posterior will not ``overfit'' in the sense that it assigns vanishing probability to values of $K$ that are much larger than the oracle complexity index $K_n$, but this is only half of the story.  We would similarly expect that the posterior would not ``underfit,'' i.e., it would assign vanishing probability to $K$ values much smaller than the oracle complexity index.  A proof of this conjecture requires a more delicate analysis, comparable to the variable selection consistency results that have been established for posterior distributions in linear regression settings.  
To our knowledge, there are no such concentration rates currently available in the literature for the marginal Gibbs posteriors of model complexity indices.  So, if one could understand the posterior asymptotics for $K$ in this problem well enough to establish this, then that would provide important insights for Gibbs model selection performance in general.

\section*{Acknowledgments}

This work is partially supported by the U.S.~National Science Foundation, DMS--1811802.

\appendix

\section{Proofs of the theorems}

\subsection{Proof of Theorem~\ref{thm:K.increasing}}
\label{S:proof2}

The proofs of Theorems~\ref{thm:K.fixed}--\ref{thm:K.increasing} are virtually the same, so here we prove the more general/interesting of the two results, namely, Theorem~\ref{thm:K.increasing}.  This proof depends on three lemmas (see below), whose proofs are postponed until Appendix~\ref{S:lemmas}.  

Start by making a minor re-write of the Gibbs posterior distribution, i.e., 
\begin{equation}
\label{eq:gibbs.proof1}
\Pi_n(A) = \frac{N_n(A)}{D_n} = \frac{\int_{\{\theta_K \in \Theta_K: \psi_{\theta_K} \in A\}} e^{-\omega t_n \{R_{n,K}(\theta_K) - R_{n,K}(\theta_K^\perp)\}} \, \tilde\Pi_K(d\theta_K)}{\int_{\Theta_K} e^{-\omega t_n \{R_{n,K}(\theta_K) - R_{n,K}(\theta_K^\perp)\}} \, \tilde\Pi_K(d\theta_K)}. 
\end{equation}
The difference between this version and the one presented above is that the exponent in both the numerator and denominator has an extra term $\omega t_n R_{n,K}(\theta_K^\perp)$.  Since these adjustments correspond to constant multiplicative factors in the respective integrands, including these does not affect the ratio.  Of course, this rewrite holds regardless of whether we are in the fixed-$K$ or increasing-$K$ cases.  

Then the strategy is to suitably lower-bound the denominator, $D_n$, and upper-bound the numerator, $N_n(A_n)$, where the set $A_n$ would appear to be different in the two theorems.  However, the two versions of $A_n$ are effectively the same.  Indeed, for the fixed-$K$ case, we have 
\begin{align}
A_n & = \{\theta_K \in \Theta_K: \|\psi_{\theta_K} - \psi_K^\perp\|_{L_2(D)} > M_n \eps_n\} \notag \\
& = \{\theta_K \in \Theta_K: \|\theta_K - \theta_K^\perp\| > M_n \eps_n\}, \label{eq:A.K}
\end{align}
with $\eps_n = t_n^{-1/2}$ and $M_n$ an arbitrary diverging sequence.  For the $K=K_n$ case, one obvious  difference is that we seek concentration around $\psi^\star$.  But recall that $\psi^\star$ is assumed to belong to the Besov space with known smoothness index $\alpha$, and the chosen $K_n = \ceiling{t_n^{1/(2\alpha+1)}}$ depends on that known $\alpha$.  Also, recall what motivated this particular choice of $K_n$, i.e., the approximation result \eqref{eq:approx}, which implies that 
\begin{align*}
\|\psi_{\theta_{K_n}} - \psi^\star\|_{L_2(D)} > M_n \eps_n & \implies \|\psi_{\theta_{K_n}} - \psi_{K_n}^\perp\|_{L_2(D)} - \|\psi_{K_n}^\perp - \psi^\star\|_{L_2(D)} \\
& \implies \|\theta_{K_n} - \theta_{K_n}^\perp\| > M_n' \eps_n, 
\end{align*}
where $\eps_n = t_n^{-\alpha/(2\alpha+1)}$ is the minimax optimal rate, and $M_n'$ is simply $M_n$ minus a constant, hence is just another sequence diverging arbitrarily slowly.  Therefore, 
\begin{equation}\label{eq:area_An_prime}
    A_n = \{\theta_{K_n}: \|\psi_{\theta_{K_n}} - \psi^\star\|_{L_2(D)} > M_n \eps\} \subseteq \{\theta_{K_n}: \|\theta_{K_n} - \theta_{K_n}^\perp\| > M_n' \eps_n\}.
\end{equation}
If the set on the right-hand side is denoted by $A_n'$, then we immediately see that it is exactly of the same form as that \eqref{eq:A.K}.  So, it suffices for us to prove that $\Pi_n(A_n')$ is vanishing.  Henceforth, we will drop the ``prime'' notation and work with 
\begin{equation}
    A_n = \{\theta_{K_n}: \|\theta_{K_n} - \theta_{K_n}^\perp\| > M_n \eps_n\}, 
\end{equation}
with $K_n$ and $\eps_n$ as defined above.  

First a bit of notation and a preliminary result concerning the properties of the projection estimator.  From the discussion in Section~\ref{SS:old}, at least in the fixed-$K$ case, recall that $\hat\theta_K$ is a consistent estimator of $\theta_K^\perp$ as $n \to \infty$ and $\Delta \to 0$.  Similar properties can be established in the case of $K=K_n$; see the proof of Lemma~\ref{lem:event.2}.  We want to highlight the event where $\hat\theta_{K_n}$ is sufficiently close to $\theta_{K_n}^\perp$, so we define 
\begin{equation}
\label{eq:event.2}
\event_n = \{ \|\hat\theta_{K_n} - \theta_{K_n}^\perp\| \leq L_n \eps_n\}.
\end{equation}
This depends on a deterministic sequence $L_n \to \infty$; the particular choice is arbitrary, and will be specified later in the proof below. The first observation is that the event $\event_n$ has probability converging to 1.  

\begin{lem}
\label{lem:event.2}
Under Condition~\ref{cond:K.increasing}, if $L_n \to \infty$, then $\prob(\event_n) \to 1$.
\end{lem}

The importance of this observation is that, on the event $E_n$, we can conveniently control the Gibbs posterior numerator and denominator.  

\begin{lem}
\label{lem:den.2}
If Condition~\ref{cond:K.increasing} holds, then $D_n$ satisfies 
\[ D_n^{-1} \, 1(\event_n) \lesssim \eps_n^{-K_n} e^{(1 + 2 L_n) \omega t_n \eps_n^2}. \]
\end{lem}

\begin{lem}
\label{lem:num.2}
If Condition~\ref{cond:K.increasing} holds, and if the $L_n$ in \eqref{eq:event.2} is chosen to satisfy $L_n \ll M_n$, where $M_n$ is the sequence in $A_n$, then  $N_n(A_n)$ satisfies 
\[ N_n(A_n) \, 1(\event_n) \lesssim (C\eps_n)^{K_n}  e^{-\rho^2\omega M_n^2t_n \eps_n^2} . \]
for any $\rho<1$ and a constant $C$ depending on $\omega$ and $\rho$
\end{lem}

Now we can proceed with the proof of Theorem~\ref{thm:K.increasing}.  Of course, $\Pi_n(A_n)$ is a probability, so it is upper bounded by 1.  Then  
\begin{align*}
\E\,\Pi_n(A_n) & = \E\{ \Pi_n(A_n) \, 1(\event_n) + \Pi_n(A_n) \, 1(\event_n^c)\} \\
& \leq \E\{ D_n^{-1} N_n(A_n) \, 1(\event_n)\} + \prob(\event_n^c).
\end{align*}
The second term in the upper bound above is $o(1)$, as $n \to \infty$, by Lemma~\ref{lem:event.2}.  By Lemmas~\ref{lem:den.2} and \ref{lem:num.2}, the first term is upper bounded by 
\[ C^{K_n} e^{-\omega t_n \eps_n^2(\rho^2 M_n^2 - 1 - 2L_n)} = \exp\{-(\rho^2 M_n^2 - 1 - 2L_n - \omega^{-1} \log C) \, \omega t_n \eps_n^2\}, \]
where the equality used the fact that $K_n = t_n \eps_n^2$.  Since $M_n$ is much larger than $L_n$, so the quantity in parentheses above is positive and, hence, the right-hand side in the above display is $o(1)$ as $n \to \infty$.  Then both terms in the upper bound for $\E\,\Pi_n(A_n)$ are vanishing, which completes the proof of Theorem~\ref{thm:K.increasing}. 

\subsection{Proof of Theorem~\ref{thm:K.unknown}}
\label{S:proof3}

In the case where a prior is introduced for the complexity index $K$, it will be helpful to rewrite the Gibbs posterior distribution as 
\begin{equation}
\label{eq:gibbs.alt}
\Pi_n(A) = \frac{N_n(A)}{D_n} = \frac{ \sum_K \tilde\pi(K) \, H_{n,K} \, N_{n, K}(A) }{ \sum_K \tilde\pi(K) \, H_{n,K} \, D_{n,K}}, 
\end{equation}
where $D_{n,K}$ and $N_{n,K}$ are as before, in particular, 
\begin{align*}
D_{n,K} & = \int_{\Theta_K} e^{-\omega t_n \{R_{n,K}(\psi_\theta) - R_{n,K}(\psi_K^\perp)\}} \, \tilde\Pi_K(d\theta) \\
N_{n,K}(A) & = \int_{A \cap \Theta_K} e^{-\omega t_n \{R_{n,K}(\psi_\theta) - R_{n,K}(\psi_K^\perp)\}} \, \tilde\Pi_K(d\theta),
\end{align*}
but now we have the additional factor 
\[ H_{n,K} = \exp\{-\omega t_n\{R_{n,K}(\psi^\perp_K)-R_{n,K_n}(\psi^\perp_{K_n})-\|\psi_{K_n}^\perp\|^2+\|\psi^\star\|^2\}\}. \]
As before, $K_n = \ceiling{ t_n^{1/(2\alpha+1)} }$, where $\alpha$ is the (unknown) smoothness index.  Note that $R_{n,K_n}(\psi_{K_n}^\perp)$, $\|\psi_{K_n}^\perp\|$, and $\|\psi^\star\|$ are constants, i.e., they do not depend on $K$ or on a generic $\theta$, so including it in the numerator and denominator does not affect the ratio.  

Next, a bit more notation.  For the arbitrary diverging sequence $M_n \to \infty$ and the target concentration rate $\eps_n = (\log t_n)^{1/2} t_n^{-\alpha/(2\alpha+1)}$ from the theorem, define
\begin{equation}
\label{eq:kn}
k_n = \max\{K: \|\psi_K^\perp - \psi^\star\| > M_n \eps_n\},
\end{equation}
the largest integer $K$ such that $\|\psi_K^\perp - \psi^\star\| > M_n \eps_n$.  Note that $k_n$ is increasing but slower than that ``ideal'' complexity $K_n = \ceiling{ t_n^{1/(2\alpha+1)} }$ from before, due to the presence of $M_n \to \infty$ in \eqref{eq:kn} and the logarithmic term in $\eps_n$.  Now define the following sequence $K_n'$, which is slightly larger than $K_n$:
\[ K_n' =\ceiling{K_n \lambda^2_n}, \quad \text{where} \quad \lambda_n = (\log t_n)^{1/2} \to \infty. \]
Finally, split the support of prior distribution $\tilde{\pi}$ into three ranges:
\begin{align*}
\K_1 & = \{K: 0 < K \leq k_n\} \\
\K_2 & = \{K: k_n < K \leq K_n' \} \\
\K_3 & = \{K: K_n' < K \leq \ceiling{t_n}\}.%
\end{align*}

Next, we define three events, which we will denote by $\event_{n,1}$, $\event_{n,2}$, and $\event_{n,3}$, corresponding to the three ranges of $K$.  The first two are 
\begin{align*}
\event_{n,1} & = \{ \|\hat\psi_{K_n} - \psi_{K_n}^\perp\| \leq L_n \lambda_n^{-1} \eps_n\} \\
\event_{n,2} & = \{ \|\hat\psi_{K_n'} - \psi_{K_n'}^\perp\| \leq L_n \eps_n\}
\end{align*}
where, as before, $L_n$ is an arbitrary diverging sequence.  The third event is 
\[ \event_{n,3} = \bigcap_{K \in \K_3} \{ \|\hat\psi_K - \psi_K^\perp\| \leq C \delta_n(K)\}, \]
where $C > 0$ is the constant specified in Condition~\ref{cond:K.unknown}.5 and 
\[ \delta_n^2(K) = t_n^{-1} K \log K. \]
Note that the target rate $\eps_n$ satisfies $\eps^2_n = K'_n/t_n$.  

To compare briefly with the corresponding setup in the proof of Theorem~\ref{thm:K.increasing}, here we are using the $L_2(D)$ norm on the basis function linear combinations, rather than the $\ell_2$ norm on the basis coefficients, but this is just for notational convenience later; the basis functions are orthonormal, so these two norms are actually the same.  The more important observation is related to the assume nestedness of the sieves.  Indeed, on event $\event_{n,1}$, we can conclude that 
\[ \|\hat\psi_K - \psi_K^\perp\| \leq L_n \lambda_n^{-1} \eps_n, \quad \text{for all $K \leq K_n $}, \]
and, similarly, on event $\event_{n,2}$, we can conclude 
\[ \|\hat\psi_K - \psi_K^\perp\| \leq L_n \eps_n, \quad \text{for all $K \leq K'_n $}. \]
That is, nestedness implies the bounds in $\event_{n,j}$ hold uniformly for all $K \in \K_j$, $j=1,2$.  The third event also implies uniform control, but this is made explicit in the intersection because it will be important in what follows that we have $K$-specific bounds.  

The essential point is that, on the event $\event_{n,1} \cap \event_{n,2} \cap \event_{n,3}$, we have very strong control on the relevant data-dependent features that contribute to the Gibbs posterior \eqref{eq:gibbs.alt}.  Lemma~\ref{lem:event.3} below says it is safe to restrict our attention to those ``nice'' data sets.  



\begin{lem}
\label{lem:event.3}
Define the event 
\begin{equation}
\label{eq:event.cap}
\event_n = \event_{n,1} \cap \event_{n,2} \cap \event_{n,3}.
\end{equation}
If the assumptions of Theorem~\ref{thm:K.unknown} hold, then $\prob(\event_n) \to 1$ as $n \to \infty$. 
\end{lem}

Now we can proceed to the details of the proof of Theorem~\ref{thm:K.unknown}.  Start with the denominator $D_n$ in \eqref{eq:gibbs.alt}.  Clearly we have 
\[ D_n > \tilde\pi(K_n) \, H_{n,K_n} \, D_{n,K_n}. \]
From the definition of $H_{n,K}$ above, when $K=K_n$, it is relatively easy to check---see the proof of Lemma~\ref{lem:H} below---that $H_{n,K_n} = \exp\{-\omega t_n \|\psi_{K_n}^\perp - \psi^\star\|^2\}$.  Since we have control on the approximation error when $K=K_n$, it follows immediately that 
\[ H_{n,K_n} \geq e^{-\omega t_n K_n^{-\alpha}} = e^{-\omega K_n}. \]
Compared to Lemma~\ref{lem:den.2}, $D_{n,K_n}$ is exactly the quantity there, and the event $\event_n$ under consideration here implies the event there.  So the same conclusion from Lemma~\ref{lem:den.2} holds here as well, i.e., 
\[ D_{n,K_n}^{-1} \, 1(\event_n) \lesssim \exp\{(1 + 2L_n) \omega \lambda_n^{-1} t_n \eps_n^2\}. \]
This looks slightly different because the $\eps_n$ here is slightly larger than that from Lemma~\ref{lem:den.2}.  Combining this with the bound on $H_{n,K_n}$ and the definition of $\tilde\pi(K_n)$, we get 
\begin{equation}
\label{eq:den.bound.3}
D_n^{-1} \, 1(\event_n) \lesssim e^{\{\beta + o(1)\} K_n \log K_n}, 
\end{equation}
where the ``$o(1)$'' term is a result of $\lambda_n \to \infty$ being of larger order than $1 + 2L_n$, since we are free to choose $L_n \ll M_n$ as small as we like.

Next, for the numerator in \eqref{eq:gibbs.alt}, recall that it depends on a particular subset
\[ A_n = \{\psi: \|\psi - \psi^\star\| > M_n \eps_n\}, \]
where $\eps_n$ is as above and $M_n$ is an arbitrary diverging sequence.  We propose to split the numerator sum over $K$ in \eqref{eq:gibbs.alt} into three pieces:
\[ N_n(A_n) = \Bigl( \sum_{K \in \K_1} + \sum_{K \in \K_2} + \sum_{K \in K_3} \Bigr) \tilde\pi(K) \, H_{n,K} \, N_{n,K}(A_n). \] 
The idea is that, when $K$ is small, the bias from using a too-simple model $\S_K$ to approximate $\psi^\star$ will be the dominant term; when $K$ is large, the model $\S_K$ is plenty adequate, so the prior that penalizes complex models will be the dominant term; and for all moderate $K$, the situation is very much like that in Theorem~\ref{thm:K.increasing} with $K=K_n$.  

When $\event_n$ in \eqref{eq:event.cap} holds, we can readily bound the various terms appearing in the Gibbs posterior.  The next result bounds the $H_{n,K}$ terms. 

\begin{lem}
\label{lem:H}
Under the assumptions of Theorem~\ref{thm:K.unknown},
\begin{equation}
\label{eq:H.bound}
H_{n,K} \, 1(\event_n) \leq 
\begin{cases}
e^{-(\omega M_n^2/2) t_n \eps_n^2} & \text{if $K \in \K_1$} \\
e^{2\omega L_nM_n \lambda_n^{-1} t_n \eps_n^2}& \text{if $K \in \K_2$} \\
e^{2C\omega \lambda_n^{-1} t_n \eps_n \delta_n(K)} & \text{if $K \in \K_3$},
\end{cases}
\end{equation}
where $L_n$ and $M_n$ are the sequences defined in $\event_n$ and $A_n$, respectively.  
\end{lem}


It remains to bound $N_{n,K}(A_n)$.  For $K \in \K_1$ and $K \in \K_3$, it suffices to use a trivial bound on the integrand above, i.e., that obtained by plugging in the risk minimizer:
\begin{equation}
\label{eq:trivial.N.bound}
N_{n,K}(A_n) \leq e^{-\omega t_n \{R_{n,K}(\hat\psi_K) - R_{n,K}(\psi_K^\perp)\}} = e^{\omega t_n \|\hat\psi_K - \psi_K^\perp\|^2}. 
\end{equation}
For $K \in \K_2$, we can proceed very much like we did in Lemma~\ref{lem:num.2}.  

\begin{lem}
\label{lem:num.3}
Under the assumptions of Theorem~\ref{thm:K.unknown}, the $K$-specific terms $N_{n,K}(A_n)$ in the numerator of \eqref{eq:gibbs.alt} satisfy
\[ N_{n,K}(A_n) \, 1(\event_n) \leq 
\begin{cases}
e^{\omega L_n^2 \lambda_n^{-2} t_n \eps_n^2 } & \text{if $K \in \K_1$} \\
e^{-\omega M_n^2 t_n\eps_n^2/2} & \text{if $K \in \K_2$} \\
e^{C^2\omega t_n \delta_n^2(K)} & \text{if $K \in \K_3$}.
\end{cases}
\] 
\end{lem}

For the remainder of the proof, we will restrict to the event $\event_n$ in \eqref{eq:event.cap} without explicitly stating so.  Start with the case $K \in \K_1$.  Here, neither the prior, $\tilde\pi(K)$, nor the numerator term, $N_{n,K}(A_n)$, is vanishing, so these will not help with the upper bound.  Fortunately, the $H_{n,K}$ term is small and dominates in the small-$K$ cases.   To see this, plug in the bounds from Lemmas~\ref{lem:H}--\ref{lem:num.3}, to get 
\[ \sum_{K \in \K_1} \tilde\pi(K) \, H_{n,K} \, N_{n,K}(A_n) \leq e^{-(\omega M_n^2 / 2) t_n \eps_n^2
+\omega L_n^2 \lambda_n^{-2} t_n \eps_n^2}. \]
Since $M_n \gg L_n$ and $\lambda_n^{-1} = o(1)$, it is clear that this upper bound is of the order $e^{-c M_n^2 t_n \eps_n^2}$ for a constant $c > 0$.  Furthermore, since $t_n \eps_n^2$ is proportional to $K_n \log K_n$, which matches the exponent from the denominator's contribution, $D_n^{-1}$, the fact that the numerator has $M_n \to \infty$ means the negative term dominates.  Therefore, the small-$K$ range's contribution to the posterior probability, $\Pi_n(A_n)$, is vanishing.    

Next, for the middle range, $K \in \K_2$, again the prior does not contribute to the bound; however, in this case, the $H_{n,K}$ term is not small and the $K$-specific numerator $N_{n,K}(A_n)$ is dominant.  As above, just plug in the bounds from Lemmas~\ref{lem:H}--\ref{lem:num.3} to get 
\[ \sum_{K \in \K_2} \tilde\pi(K) \, H_{n,K} \, N_{n,K}(A_n) \leq e^{-\omega M_n^2 t_n \eps_n^2/2
+2\omega L_nM_n \lambda_n^{-1} t_n \eps_n^2} \leq
e^{-c \omega M_n^2 t_n \eps_n^2},\]
for another constant $c > 0$. Just like in the previous discussion, the moderate-$K$ range's contribution to the posterior probability is also vanishing.

Finally, for $K \in \K_3$, neither $H_{n,K}$ nor $N_{n,K}(A_n)$ are small, but since the prior is decaying sufficiently fast, we can get a vanishing bound.  That is,
\begin{align}
\sum_{K \in \K_3} \tilde\pi(K) \, H_{n,K} \, N_{n,K}(A_n) & \leq
    \sum_{K \in \K_3}
    e^{-\beta K \log K} e^{C^2\omega t_n \delta_n^2(K)} e^{2C\omega \lambda_n^{-1} t_n \eps_n \delta_n(K)} \notag \\
& = \sum_{K \in \K_3} e^{-\{\beta - C^2 \omega + o(1)\} K \log K}.  \label{eq:sum_in_range3}
\end{align}
The last equality follows because 
\[ \lambda_n^{-1} t_n \eps_n \delta_n(K) \propto t_n^{\frac12 - \frac{\alpha}{2\alpha+1}} (K \log K)^{1/2} =
\Bigl( \frac{K_n}{K \log K} \Bigr)^{1/2} (K \log K), \]
and the ratio is $o(1)$ since $K \gg K_n$ in $\K_3$.  Provided that $\beta > C^2 \omega + o(1)$, as in \eqref{eq:beta}, the sum in \eqref{eq:sum_in_range3} can be bounded by the tail of a geometric series, i.e., 
\[ \sum_{K \in \K_3} e^{-\{\beta - C^2 \omega + o(1)\} K \log K} \leq \sum_{K > K_n'} (K_n')^{-\{\beta - C^2\omega + o(1)\} K} = e^{-\{\beta - C^2\omega + o(1)\} K_n' \log K_n'}. \]
The contribution from the denominator looks similar, but since $K_n' \gg K_n$, it follows that quantity above for the numerator is the dominant term.  

Putting everything together, since the contribution to $\Pi_n(A_n)$ from all three $K$ ranges is vanishing, we can apply the same logic used in the very last step of the proof of Theorem~\ref{thm:K.increasing} to show that $\E\{\Pi_n(A_n)\} \to 0$. This completes the proof of Theorem~\ref{thm:K.unknown}.

\subsection{Proof of Theorem~\ref{thm:no.overfit}}
\label{proof:no.overfit}

In general, the marginal Gibbs posterior for $K$ is given by 
\[ \pi_n(K) = \frac{\tilde\pi(K) H_{n,K} N_{n,K}(L_2)}{D_n}, \]
where $H_{n,K}$, $N_{n,K}$, and $D_n$ are defined as in the proof of Theorem~\ref{thm:K.unknown}.  The range of $K$ that we are concerned with is contained in $\K_3$, so the event $\event_n$ which, according to Lemma~\ref{lem:event.3}, has probability approaching 1, gives us the control we need on the data-dependent terms.  In particular, on $\event_n$, we have the following:
\begin{itemize}
\item a lower bound on $D_n$ of the form
\[ D_n > e^{-\{\beta + o(1)\} K_n \log K_n}; \]
\item from Lemma~\ref{lem:H}, an upper bound on $H_{n,K}$ of the form 
\[ H_{n,K} \leq e^{2C \omega \lambda_n^{-1} t_n \eps_n \delta_n(K)} \leq e^{2C \omega \lambda_n^{-1} K \log K}; \]
\item and, by Lemma~\ref{lem:num.3}, an upper bound on $N_{n,K}(L_2)$ of the form
\[ N_{n,K}(L_2) \leq e^{\omega t_n \|\hat\psi_K - \psi_K^\perp\|^2} \leq e^{C^2 \omega K \log K}. \]
\end{itemize}
Putting everything together, on event $\event_n$, 
\begin{align*}
\pi_n(K) & \leq \exp[-\beta K \log K + \{\beta + o(1)\} K_n \log K_n + (C^2 + 2C \lambda_n^{-1}) \omega K \log K] \\
& = \exp[-\{ \beta - C^2 \omega + o(1)\} K \log K + \{\beta + o(1)\} K_n \log K_n].
\end{align*}
Let $\zeta = \beta - C^2 \omega + o(1)$, which is positive by the conditions we imposed on $(\beta,\omega)$.  Then 
\begin{align*}
\Pi_n(\{K: K > \tau K_n\}) & = \sum_{K > \tau K_n} \pi_n(K) \\
& \leq e^{\{\beta + o(1)\} K_n \log K_n} \sum_{K > \tau K_n} e^{-\zeta K \log K} \\
& \leq e^{\{\beta + o(1)\} K_n \log K_n} \sum_{K > \tau K_n} (\tau K_n)^{-\zeta K} \\
& \lesssim e^{\{\beta + o(1)\} K_n \log K_n} e^{-\zeta(\tau K_n + 1) \log(\tau K_n)}. 
\end{align*}
For this to be vanishing, we need $\tau \zeta \equiv \tau(\beta - C^2 \omega) > \beta$ or, equivalently, 
\[ \beta > (\tau -1)^{-1} \tau C^2 \omega. \]
This is precisely the condition \eqref{eq:beta.new}, so under the stated assumptions, the conclusion of Theorem~\ref{thm:no.overfit} follows.

\section{Proofs of the lemmas}
\label{S:lemmas}

\subsection{Proof of Lemma \ref{lem:event.2}}

Following the proof of Lemma~3.2 in \citet{figueroa2009nonparametric} and applying it with the improved small-time bound in Proposition~2.1 in \citet{MR2787609}, we find that 
\begin{align}
\E \|\hat\theta_K - \theta_K^\perp\|^2
& =\E\|\hat\theta_K - \E\hat\theta_K \|^2+\|\E\hat\theta_K-\theta_K^\perp\|^2 \notag \\
& \lesssim K t_n^{-1} \{ 1 + F_1^2(K) \Delta^2 t_n + F_2(K) \Delta\}, \label{eq:jose.mse}
\end{align}
where $F_1(K)$ and $F_2(K)$ are those features $\S_K$ as described in Section~\ref{SS:setup}. Under the stated assumptions on $\Delta$, whether we are in the fixed-$K$ or increasing-$K$ case, the right-hand side of \eqref{eq:jose.mse} is $\lesssim K t_n^{-1}$.  Now apply Markov's inequality to get 
\[ \prob(\|\hat\theta_K - \theta_K^\perp\| > L_n \eps_n) \lesssim \frac{K}{L_n^2 t_n \eps_n^2}. \]
In the case where $K$ is a constant and $\eps_n = t_n^{-1/2}$, clearly the upper bound vanishes since $L_n \to \infty$.  For the case when $K=K_n$, recall that $K_n$ and $\eps_n$ are connected through the smoothness index, $\alpha$, of the L\'evy density $\psi^\star$, i.e., 
\[ \eps_n = t_n^{-\alpha/(2\alpha+1)} \quad \text{and} \quad K_n = \ceiling{t_n^{1/(2\alpha+1)}}. \]
Then, in this case, it is easy to confirm that 
\[ \frac{K_n}{t_n \eps_n^2} = O(1), \quad n \to \infty, \]
so, again, the upper bound vanishes since $L_n \to \infty$.




\subsection{Proof of Lemma~\ref{lem:den.2}}

For simplicity, we drop the ``$K$'' in the notation momentarily.  By direct calculation, 
\begin{align*}
D_n & = \int e^{-\omega t_n\{ R_n(\theta) - R_n(\theta^\perp)\}} \, \tilde\Pi(d\theta) \\
& \geq \int_{\|\theta-\theta^\perp\| \leq \eps_n} e^{-\omega t_n \{ R_n(\theta) - R_n(\theta^\perp)\}} \, \tilde\Pi(d\theta) \\
& \geq \tilde\Pi(\{\theta: \|\theta-\theta^\perp\| \leq \eps_n\}) \, \inf_{\|\theta-\theta^\perp\| \leq \eps_n} e^{-\omega t_n \{ R_n(\theta) - R_n(\theta^\perp)\}}. 
\end{align*}
First, we bound the prior-dependent term.  For the fixed-$K$ case, it is easy to check that 
\[ \tilde\Pi_K(\{\theta: \|\theta-\theta^\perp\| \leq \eps_n\}) \gtrsim \eps_n^K. \]
For the $K=K_n$ case, start by replacing the $\ell_2$ ball around $\theta^\perp$ with the $\ell_\infty$ ball, so that 
\[ \tilde\Pi_{K_n}(\{\theta: \|\theta-\theta^\perp\| \leq \eps_n\}) \geq \tilde\Pi_K(\{\theta: \max_k |\theta_{K_nk}-\theta_{K_nk}^\perp| \leq K_n^{1/2} \eps_n\}). \]
Since the prior treats the components independently, the right-hand side can be written as a product:
\[ \tilde\Pi_{K_n}(\{\theta: \max_k |\theta_{K_nk}-\theta_{K_nk}^\perp| \leq K_n^{1/2} \eps_n\}) = \prod_{k=1}^{K_n} \tilde\Pi_{K_nk}(\{\theta_{K_nk}: |\theta_{K_nk} - \theta_{Kk}^\perp| \leq K_n^{1/2} \eps_n\}). \]
Each term in the product on the right-hand side can be lower-bounded by a constant times $K_n^{1/2} \eps_n$.  Since $\eps_n \sim K_n^{-\alpha}$ and $\alpha > \frac12$, the intervals 
\[ |\theta_{K_nk} - \theta_{K_nk}^\perp| \leq K_n^{1/2} \eps_n, \quad k=1,\ldots,K_n, \]
are bounded, which implies that the marginal density (mass) function is bounded away from 0 on those intervals.  Therefore, there exists a constant $J \in (0,1)$ such that 
\[ \tilde\Pi_{K_n}(\{\theta: \|\theta-\theta^\perp\| \leq \eps_n\}) \gtrsim J^{K_n} \geq \eps_n^{K_n}. \]

Second, the empirical risk difference can be simplified as 
\[ R_n(\theta) - R_n(\theta^\perp) = \|\theta-\theta^\perp\|^2 - 2 \langle \theta-\theta^\perp, \hat\theta - \theta^\perp \rangle. \]
By Cauchy--Schwartz, this is upper bounded by 
\[ \|\theta-\theta^\perp\|^2 + 2 \|\theta-\theta^\perp\| \, \|\hat\theta - \theta^\perp\|. \]
On the event $\event_n$, $\|\hat\theta-\theta^\perp\|$ is upper bounded by $L_n \eps_n$, so 
\[ e^{-\omega t_n \{R_n(\theta) - R_n(\theta^\perp)\}} \geq e^{-\omega t_n (\|\theta - \theta^\perp\|^2 + 2L_n\eps_n \|\theta-\theta^\perp\|)}, \quad \text{on $\event_n$}. \]
The right-hand side is a decreasing function of $\|\theta-\theta^\perp\|$ so, 
\[ \inf_{\|\theta-\theta^\perp\| \leq \eps_n} e^{-\omega t_n\{ R_n(\theta) - R_n(\theta^\perp)\}} \geq e^{-(1+2L_n) t_n \eps_n^2}. \]
Combining the two lower bounds established above proves that $D_n^{-1}$ is suitably upper bounded on the event $\event_n$ and, hence, the claim.

\subsection{Proof of Lemma~\ref{lem:num.2}}

As above, for simplicity, we drop the ``$K$'' in the notation. Recall that $N_n(A)$ is an integral of the exponentiated empirical risk difference 
over $A$. As discussed in Section \ref{S:proof2}, it suffices to derive the upper bound for $N_n(A_n)$ with an $A_n$ given in \eqref{eq:area_An_prime}. 
Alternatively, $A_n$ can be re-expressed as a union of disjoint shells, $\bigcup_{J=1}^\infty Q_{n,J}$, where 
\[ Q_{n,J} = \{\theta: J M_n \eps_n < \|\theta-\theta^\perp\| \leq (J+1) M_n \eps_n\}, \quad J \geq 1. \]   

Then the Gibbs posterior numerator is 
\begin{align*}
N_n(A_n) & = \sum_{J=1}^\infty \int_{Q_{n,J}} e^{-\omega t_n \{R_n(\theta) - R_n(\theta^\perp)\}} \, \tilde\Pi(d\theta) \\
& \leq \sum_{J=1}^\infty \tilde\Pi(Q_{n,J}) \times \sup_{\theta \in Q_{n,J}} e^{-\omega t_n \{R_n(\theta) - R_n(\theta^\perp)\}}. 
\end{align*}
The empirical risk function is continuous and convex, with unique minimizer at $\hat\theta$.  On the event $\event_n$, $\hat\theta$ satisfies $\|\hat\theta-\theta^\perp\| \ll M_n\eps_n$, so it is not contained in any of the shells.  Therefore, on each shell, the supremum is attained on the boundary where $|\theta-\theta^\perp\| = J M_n \eps_n$, i.e.,  
\[ \sup_{\theta \in Q_{n,J}} e^{-\omega t_n \{R_n(\theta) - R_n(\theta^\perp)\}} = \sup_{\theta: \|\theta-\theta^\perp\| = J\rho'M_n \eps_n} e^{-\omega t_n \{R_n(\theta) - R_n(\theta^\perp)\}}. \]
Consequently, applying Cauchy--Schwartz to the empirical risk difference, we get that 
\[ \sup_{\theta \in Q_{n,J}} e^{-\omega t_n \{R_n(\theta) - R_n(\theta^\perp)\}} \leq e^{-\omega t_n M_n^2 \eps_n^2 J^2 (1 - 2L_n / J M_n)} \leq e^{-\rho\omega t_n M_n^2 \eps_n^2 J^2}, \]
where the last inequality follows from facts that $J \geq 1$ and $1-2L_n/M_n \geq\rho$ when $n$ is enough large. Since the prior has a bounded density, the shell probabilities satisfy 
\[ \tilde\Pi(Q_{n,J}) \lesssim (M_n \eps_n)^{K_n} (J+1)^{K_n} \leq (M_n \eps_n)^{K_n}(2J)^{K_n}, \quad J=1,2,\ldots. \]
Plugging these two bounds into the above expression, we have 
\[ N_n(A_n) \, 1(\event_n) \lesssim (M_n \eps_n)^{K_n} \sum_{J=1}^\infty (2J)^{K_n} e^{-\rho\omega t_n M_n^2 \eps_n^2 J^2}. \]
Recall that $K_n$ and $t_n \eps_n^2$ are equivalent.  The sum is asymptotically equivalent to the integral below which, after some appropriate change-of-variables, can be bounded as 
\[ \int_1^\infty (2u)^K e^{-\rho\omega M_n^2 K_n u^2} \, du \lesssim 2^{K_n-1} (\rho\omega M_n^2 K_n)^{-(K_n+1)/2} \,  \Gamma(\tfrac{K_n+1}{2}, \rho\omega M_n^2 K_n), \]
where $\Gamma(s, x) = \int_x^\infty y^{s-1} e^{-y} \, dy$ is the incomplete gamma function.  Since $s \mapsto \Gamma(s,x)$ is increasing, we can assume in what follows that $K_n$ is an odd number. 


If $K_n$ is odd, then $(K_n+1)/2$ is an integer, and the incomplete gamma function has the following expansion:
\begin{equation}
\label{eq:inga_sum}
\Bigl(\frac{K_n-1}{2} \Bigr)! \times e^{-(K_n\tau_n+\tau_n)} \sum_{r=0}^{(K_n-1)/2} \frac{(K_n\tau_n+\tau_n)^r}{r!},
\end{equation}
where $\tau_n=\rho\omega M_n^2 K_n(K_n+1)^{-1}$ is a  sequence $\to \infty$ at the same speed as $M^2_n$. In \eqref{eq:inga_sum}, the expression to the right of ``$\times$'' is the probability of a Poisson random variable, with rate $K_n\tau_n+\tau_n$, being smaller than $(K_n-1)/2$. An upper bound on this cumulative probability, according to \citet{short2013improved}, is given as,
\[ \Phi\Bigl[\text{sign}\bigl\{\tfrac12(K_n+1)(1-2\tau_n)\bigr\} \cdot \bigl\{ (K_n+1)(2\tau_n-\log2\tau_n-1)
\bigr\}^{1/2} \Bigr], \]
where $\Phi$ is the standard normal distribution function. Since $\tau_n$ is a diverging sequence, the above sign function has a negative value, and the lower tail probability of a standard normal distribution is bounded by $e^{-\rho(K_n+1)\tau_n}=e^{-\rho^2\omega M_n^2K_n}$, when $n$ is enough large. 
Combining the above results, we obtain
\begin{align*}
N_n(A_n) \, 1(\event_n) \lesssim & \eps_n^{K_n} (M_n)^{-1}
(\rho\omega/2)^{-(K_n+1)/2}
 e^{-\rho^2\omega M_n^2K_n} \frac{(K_n/2)!}{(K_n/2)^{K_n/2}}\\
\lesssim & \eps_n^{K_n} 
(\rho\omega/2)^{-K_n/2}
 e^{-\rho^2\omega M_n^2K_n},
\end{align*}
then Lemma \ref{lem:num.2} holds with constant $C=(\rho\omega/2)^{-1/2}$.


\subsection{Proof of Lemma~\ref{lem:event.3}}

With only minor modifications, we can follow the proof of Lemma~\ref{lem:den.2} to show that $\prob(\event_{n,1}) \to 1$ and $\prob(\event_{n,2}) \to 1$.  So here it remains to show that $\prob(\event_{n,3}) \to 1$.  

Towards this, first we show that it suffices for us to show the following event 
\begin{align}\label{eq:event_E3_suffice}
    \event_{n,3}' = \bigcap_{K \in \K_3} \underbrace{ \{ \|\hat\psi_{(K-1):K} - \E \psi_{(K-1):K}^\perp\| \leq C' K^{-1/2}\delta_n(K)\} }_{\event_{n,3}(K)}.
\end{align}
has probability increasing to $1$, where $\hat\psi_{K':K} := \hat\psi_K - \hat\psi_{K'}$, and $C'$ is an arbitrary constant smaller than $C$. 
We note that $\event_{n,3}'$ is the intersection of a sequence of events $\event_{n,3}(K)$, for $K\in\K_3$, and each $\event_{n,3}(K)$ places a $K$-depending upper bound on the variability introduced by basis function $f_K$. 
To see that (\ref{eq:event_E3_suffice}) is a sufficient claim, we start with the following seemingly trivial decomposition, 
\[ \hat\psi_K = \hat\psi_{K_n'} + \hat\psi_{K_n':K}, \quad K > K_n'. \] 
This will be helpful because, since the sieves are nested, the function set $\S_{K_n':K}$ that contains $\hat\psi_{K_n':K}$ is much simpler than $\S_K$, when $K > K_n'$.  
By the triangle inequality, 
\[ \|\hat\psi_K - \psi_K^\perp\| \leq \|\hat\psi_{K_n':K} - \E \hat\psi_{K_n':K}\| + \|\hat\psi_{K_n'} - \E \hat\psi_{K_n'}\| + \|\E \hat\psi_K - \psi_K^\perp\|. \]
Given that events $\event_{n,2}$ and $\event_{n,3}$ hold with probabilities converging to $1$, to handle the three terms on the right-hand side, we note (a) the deterministic last term has an upper bound which is of smaller order than $\delta_n(K)$, and (b) with a ``large probability'' accompanying  $\event_{n,2}$ and $\event_{n,3}$, the other two terms have upper bounds which is of smaller/same order than/as $\delta_n(K)$.  More specifically, 
\begin{itemize}
\item From \eqref{eq:jose.mse}, we have that
\[ \|\E\hat\psi_K - \psi_K^\perp\| \lesssim K^{1/2} F_1(K) \Delta, \]
and Condition~\ref{cond:K.unknown}.2 implies that the right-hand side is $\ll \delta_n(K)$.
\item On the event $\event_{n,2}$, which has probability converging to 1, we have 
\[ \|\hat\psi_{K_n'} - \E \hat\psi_{K_n'}\| \leq L_n  \eps_n. \]
Since we are free to choose $L_n$ as small as we would like, by choosing $L_n \ll \lambda_n$, we easily see the second term on the right-hand side above is also $\ll \delta_n(K)$. 
\item On the event $\event_{n,3}'$, we have   
\begin{align}
\nonumber
\|\hat\psi_{K_n':K} - \E \hat\psi_{K_n':K}\|^2 &=\sum_{k=K_n'+1}^{K}\|\hat\psi_{(k-1):k} - \E \hat\psi_{(k-1):k}\|^2 \\\nonumber
& \leq C^{\prime 2} \sum_{k=K_n'+1}^{K}k^{-1}\delta^2_n(k) \\\nonumber
& \leq C^{\prime 2} t_n^{-1} \sum_{k=K_n'+1}^{K} \log k \\
\label{eq:lemma4_bound1}
& \leq C^{\prime 2} t_n^{-1} K \log K.
\end{align}
So the first term on the right-hand side above has an upper bound which is of the same order as $\delta_n(K)$ up to a slightly small constant $C'<C$. 
\end{itemize}

Recalling the decomposition of $\event_{n,3}'$, so to show $\prob(\event_{n,3}') \to 1$, it suffices to show that 
\begin{equation}
\label{eq:event.goal}
\sum_{K \in \K_3} \prob\{ \event_{n,3}(K)^c \} = o(1), \quad n \to \infty. 
\end{equation}
We will do so below by bounding the individual $\prob\{ \event_{n,3}(K)^c \}$ terms. Let 
\[ U_K:=\|\hat\psi_{(K-1):K} - \E \hat\psi_{(K-1):K}\| = \frac{1}{t_n} \Bigl| \sum_{i=1}^n (f_K - \E f_K)(Y_i) \Bigr|, \]
Since the basis functions are bounded by $B$, $t_n U_K$ is a sum of  iid random variables, each bounded by $2B$.  Moreover, for the variance of $t_n U_K = \sum_{i=1}^n (f_K - \E f_K)(Y_i)$, we can follow the proof of Proposition~3.4 in \citet{figueroa2009nonparametric} to get
\[ \V(t_n U_K) \leq t_n \Bigl\{ \int_D f_K^2(x) \, \psi^\star(x) \, dx + F_2(K) \Delta \Bigr\} \leq t_n \{\xi + F_2(K) \Delta\}, \]
where $\xi = \sup_{x \in D} \psi^\star(x)$ is finite by assumption. With this information, we can apply Bernstein's inequality to get 
\begin{align*}
\prob\{ |U_K| > C' K^{-1/2} \delta_n(K) \} & = \prob\Bigl\{ \Bigl| \sum_{i=1}^n (f_K - \E f_K)(Y_i) \Bigr| > C' (t_n \log K)^{1/2} \Bigr\} \\
& \leq 2\exp\Bigl[ -\frac{C^{\prime 2} \log K}{2 \{\xi + F_2(K) \Delta\} + \frac23 B C' t_n^{-1/2} (\log K)^{1/2}} \Bigr].
\end{align*}
The denominator of the above exponent is $2\xi + o(1)$, so the entire right-hand side can be written as $2K^{-\gamma}$, where 
\[ \gamma = \frac{C^{\prime 2}}{2\xi + o(1)}. \]
So long as we take $C^{\prime 2} > 2 \xi + o(1)$, which is feasible given the constraint $C' < C$ and the definition of $C$ in Condition~\ref{cond:K.unknown}.4, we have that $\gamma > 1$.  Therefore, 
\[ 2 \sum_{K \in \K_3} K^{-\gamma} \lesssim \sum_{K=K_n'+1}^\infty K^{-\gamma } = o(1), \quad n \to \infty, \]
proving \eqref{eq:event.goal} and, hence, the lemma.

\subsection{Proof of Lemma~\ref{lem:H}}

Start by defining the functions $Z_{n,K}=\hat\psi_K-\psi^\perp_K$, depending on data.  With this notation, we can write 
\begin{align*}
R_{n,K}(\psi^\perp_K) & = -\|\psi^\perp_K\|^2-2\langle\psi^\perp_K,Z_{n,K}\rangle \\
-\|\psi_{K_n}\|^2-R_{n,K_n}(\psi^\perp_{K_n}) & = 2\langle\psi^\perp_{K_n},Z_{n,K_n}\rangle. 
\end{align*}
With these, and an application of the Pythagorean theorem, we can re-express $H_{n,K}$ as 
\begin{align}
H_{n,K}=& \exp[-\omega t_n\{R_{n,K}(\psi^\perp_K)-R_{n,K_n}(\psi^\perp_{K_n})-\|\psi_{K_n}\|^2+\|\psi^\star\|^2)\}] \notag \\
=& \exp[-\omega t_n\|\psi_K^\perp - \psi^\star\|^2+2\omega t_n \{ \langle Z_{n,K}, \psi_K^\perp \rangle - \langle Z_{n,K_n}, \psi_{K_n}^\perp \rangle \}] \label{eq:Hn_decmp}
\end{align}
For the first term in the exponent of \eqref{eq:Hn_decmp}, we always have $\|\psi_K^\perp - \psi^\star\| \geq 0$; however, when $K \in \K_1$, i.e., when $K < k_n$, we know that $\|\psi_K^\perp - \psi^\star\|\geq M_n \eps_n$, which will be useful. 

For the second term in the exponent, we consider reflecting the nested structure of the sieves on the inner products,
and then applying Cauchy--Schwartz inequality:
\begin{align*} 
|\langle Z_{n,K}, \psi_K^\perp \rangle - \langle Z_{n,K_n}, \psi_{K_n}^\perp \rangle| & = |\langle Z_{n,K_n} - Z_{n,K}, \psi_{K_n}^\perp - \psi_K^\perp \rangle|\\
& \leq \|Z_{n,K_n} - Z_{n, K}\| \, \|\psi_{K_n}^\perp - \psi_K^\perp\|.
\end{align*}
Further, $\|\psi_{K_n}^\perp - \psi_K^\perp\|$ has an upper bound depending on the distance between $\psi^\star$ and its
projection onto $\S_{K\wedge K_n}$. That is
\begin{equation}
\label{eq:distance_two_proj}
    \|\psi_{K_n}^\perp - \psi_K^\perp\|^2 = \|\psi_{K\wedge K_n}^\perp - \psi^\star\|^2 - \|\psi_{K\vee K_n}^\perp - \psi^\star\|^2\leq\|\psi_{K\wedge K_n}^\perp - \psi^\star\|^2.
\end{equation}
The upper bound for $\|Z_{n,K_n} - Z_{n, K}\|$, on the other hand, is given by Lemma~\ref{lem:event.3}
\begin{equation}
\label{eq:distance_two_Z_nK}
    \|Z_{n,K_n} - Z_{n, K}\|\leq
\begin{cases}
  L_n \lambda_n^{-1} \eps_n, &\text{ if }K\leq K_n\\
  L_n \eps_n, &\text{ if } K_n<K\leq K_n'\\
  C\delta_n(K), &\text{ if }  K_n' <K\leq\ceiling{t_n}
\end{cases}
\end{equation}

To put everything together, we consider the three ranges of $K$ separately.  In all of what follows, we are implicitly restricting to the event $E_n$.  
\begin{itemize}
\item If $K \in \K_1$, then we know that $K \leq k_n < K_n$, where $k_n$ is defined in \eqref{eq:kn}.   Then 
\[ H_{n,K} \leq e^{-\omega t_n \|\psi_K^\perp - \psi^\star\|^2 + 2\omega t_n L_n \lambda_n^{-1}\eps_n \|\psi_K^\star - \psi^\star\|}. \]
Focusing specifically on the exponent, it can be rewritten as 
\[ \Bigl( \frac{2L_n \eps_n}{\lambda_n \|\psi_K^\perp - \psi^\star\|} - 1 \Bigr) \omega t_n \|\psi_K^\perp - \psi^\star\|^2 . \]
For $K$ in this range, we have $\|\psi_K^\perp - \psi^\star\| > M_n \eps_n$, so the above display is less than 
\[ \Bigl( \frac{2L_n}{M_n \lambda_n} - 1 \Bigr) \omega t_n \|\psi_K^\perp - \psi^\star\|^2 , \]
clearly, the term in parentheses above would be less than, say, $-1/2$, for all large $n$.  With the factor being negative, we can use the lower bound on the approximation error again to get that 
\[ H_{n,K} \leq e^{-(\omega M_n^2/2) t_n \eps_n^2}. \]

\item If $K \in \K_2$, since $\|\psi^\perp_K-\psi^\star\|$ is non-negative, by plugging in the upper bounds (\ref{eq:distance_two_proj}) and (\ref{eq:distance_two_Z_nK}), we get
\begin{align*}
H_{n,K} & \leq 
\begin{cases}
  e^{2\omega L_nM_n\lambda_n^{-1} t_n\eps_n^2}, &\text{ for }k_n<K\leq K_n\\
  e^{2\omega L_n \lambda_n^{-1} t_n \eps_n^2}, &\text{ for }K_n<K\leq K_n'
\end{cases} \\
& \leq e^{2 \omega L_n M_n \lambda_n^{-1} t_n \eps_n^2}.
\end{align*}

\item If $K \in \K_3$, still because  $\|\psi^\perp_K-\psi^\star\|$ is non-negative, we have
\[ H_{n,K} \leq e^{2C\omega \lambda_n^{-1} t_n \eps_n \delta_n(K)}, \]

\end{itemize}

\subsection{Proof of Lemma~\ref{lem:num.3}}

The cases $K \in \K_1$ and $K \in \K_3$ are straightforward; just plug in the bounds determined by the event $\event_n$ into \eqref{eq:trivial.N.bound}. For the $K \in \K_2$ case, we proceed like in the proof of Lemma~\ref{lem:num.2}.  In fact, the argument here is simpler because the attainable rate already has a logarithmic term in it, so we do not need the decomposition of $A_n$ into shells as we did in the proof of Lemma~\ref{lem:num.2}.  So, since the empirical risk function is convex, and since the event $E_{n,K}$ ensures that the minimizer, $\hat\psi_K$, is in the interior of $A_n^c$, it can be shown, as we did previously, that 
\[ \sup_{\psi \in A_n} e^{-\omega t_n\{R_{n,K}(\psi) - R_{n,K}(\psi_K^\perp)\}} = \sup_{\psi \in \partial A_n} e^{-\omega t_n\{R_{n,K}(\psi) - R_{n,K}(\psi_K^\perp)\}} \leq e^{-\omega t_n M_n^2 \eps_n^2/2}. \]

\bibliographystyle{apalike}
\bibliography{Template}

\begin{thebibliography}{}

\bibitem[A\"it-Sahalia and Jacod, 2014]{sahalia.jacod.book}
A\"it-Sahalia, Y. and Jacod, J. (2014).
\newblock {\em High-Frequency {F}inancial {E}conometrics}.
\newblock Princeton {U}niversity {P}ress.

\bibitem[Applebaum, 2009]{MR2512800}
Applebaum, D. (2009).
\newblock {\em L\'{e}vy Processes and Stochastic Calculus}, volume 116 of {\em
  Cambridge Studies in Advanced Mathematics}.
\newblock Cambridge University Press, Cambridge, second edition.

\bibitem[Arbel et~al., 2013]{MR3091697}
Arbel, J., Gayraud, G., and Rousseau, J. (2013).
\newblock Bayesian optimal adaptive estimation using a sieve prior.
\newblock {\em Scandinavian Journal of Statistics}, 40(3):549--570.

\bibitem[Barndorff-Nielsen et~al., 2001]{obn.levy.collection.2001}
Barndorff-Nielsen, O.~E., Mikosch, T., and Resnick, S.~I., editors (2001).
\newblock {\em L\'evy {P}rocesses: {T}heory and {A}pplications}. Birkh\"auser
  Basel.

\bibitem[Barron et~al., 1999]{MR1679028}
Barron, A., Birg\'{e}, L., and Massart, P. (1999).
\newblock Risk bounds for model selection via penalization.
\newblock {\em Probability Theory and Related Fields}, 113(3):301--413.

\bibitem[Belomestny et~al., 2019]{belomestny2018nonparametric}
Belomestny, D., Gugushvili, S., Schauer, M., and Spreij, P. (2019).
\newblock Nonparametric {B}ayesian inference for gamma-type {L}\'{e}vy
  subordinators.
\newblock {\em Communications in Mathematical Sciences}, 17(3):781--816.

\bibitem[Bertoin, 1996]{MR1406564}
Bertoin, J. (1996).
\newblock {\em L\'{e}vy Processes}, volume 121 of {\em Cambridge Tracts in
  Mathematics}.
\newblock Cambridge University Press, Cambridge.

\bibitem[Birg\'{e} and Massart, 1997]{MR1462939}
Birg\'{e}, L. and Massart, P. (1997).
\newblock From model selection to adaptive estimation.
\newblock In {\em Festschrift for {L}ucien {L}e {C}am}, pages 55--87. Springer,
  New York.

\bibitem[Bissiri et~al., 2016]{MR3557191}
Bissiri, P.~G., Holmes, C.~C., and Walker, S.~G. (2016).
\newblock A general framework for updating belief distributions.
\newblock {\em Journal of the Royal Statistical Society. Series B. Statistical
  Methodology}, 78(5):1103--1130.

\bibitem[Black and Scholes, 1973]{MR3363443}
Black, F. and Scholes, M. (1973).
\newblock The pricing of options and corporate liabilities.
\newblock {\em Journal of Political Economy}, 81(3):637--654.

\bibitem[Comte and Genon-Catalot, 2009]{MR2565560}
Comte, F. and Genon-Catalot, V. (2009).
\newblock Nonparametric estimation for pure jump {L}\'{e}vy processes based on
  high frequency data.
\newblock {\em Stochastic Processes and their Applications},
  119(12):4088--4123.

\bibitem[Cont and Tankov, 2004]{MR2042661}
Cont, R. and Tankov, P. (2004).
\newblock {\em Financial Modelling with Jump Processes}.
\newblock Chapman \& Hall/CRC Financial Mathematics Series. Chapman \&
  Hall/CRC, Boca Raton, FL.

\bibitem[Figueroa-L\'{o}pez, 2009]{figueroa2009nonparametric}
Figueroa-L\'{o}pez, J.~E. (2009).
\newblock Nonparametric estimation for {L}\'evy models based on
  discrete-sampling.
\newblock {\em IMS Lecture {N}otes-{M}onograph {S}eries}, 57:117--146.

\bibitem[Figueroa-L\'{o}pez, 2011]{MR2787609}
Figueroa-L\'{o}pez, J.~E. (2011).
\newblock Sieve-based confidence intervals and bands for {L}\'{e}vy densities.
\newblock {\em Bernoulli}, 17(2):643--670.

\bibitem[Grenander, 1981]{MR599175}
Grenander, U. (1981).
\newblock {\em Abstract {I}nference}.
\newblock Wiley Series in Probability and Mathematical Statistics. John Wiley
  \& Sons, Inc., New York.

\bibitem[Gr\"{u}nwald and van Ommen, 2017]{MR3724979}
Gr\"{u}nwald, P. and van Ommen, T. (2017).
\newblock Inconsistency of {B}ayesian inference for misspecified linear models,
  and a proposal for repairing it.
\newblock {\em Bayesian Analysis}, 12(4):1069--1103.

\bibitem[Gr\"{u}nwald and Mehta, 2020]{MR4095335}
Gr\"{u}nwald, P.~D. and Mehta, N.~A. (2020).
\newblock Fast rates for general unbounded loss functions: from {ERM} to
  generalized {B}ayes.
\newblock {\em Journal of Machine Learning Research (JMLR)}, 21:Paper No. 56,
  80.

\bibitem[Gugushvili et~al., 2015]{MR3356921}
Gugushvili, S., van~der Meulen, F., and Spreij, P. (2015).
\newblock Nonparametric {B}ayesian inference for multidimensional compound
  {P}oisson processes.
\newblock {\em Modern Stochastics. Theory and Applications}, 2(1):1--15.

\bibitem[Gugushvili et~al., 2018]{MR3769832}
Gugushvili, S., van~der Meulen, F., and Spreij, P. (2018).
\newblock A non-parametric {B}ayesian approach to decompounding from high
  frequency data.
\newblock {\em Statistical Inference for Stochastic Processes}, 21(1):53--79.

\bibitem[It\^{o}, 1942]{MR14629}
It\^{o}, K. (1942).
\newblock On stochastic processes. {I}. ({I}nfinitely divisible laws of
  probability).
\newblock {\em Japanese Journal of Mathematics}, 18:261--301.

\bibitem[L{\'e}vy, 1934]{levy1934integrales}
L{\'e}vy, P. (1934).
\newblock Sur les int{\'e}grales dont les {\'e}l{\'e}ments sont des variables
  al{\'e}atoires ind{\'e}pendantes.
\newblock {\em Annali della {S}cuola {N}ormale {S}uperiore di {P}isa-{C}lasse
  di {S}cienze}, 3(3-4):337--366.

\bibitem[Lyddon et~al., 2019]{MR3949315}
Lyddon, S.~P., Holmes, C.~C., and Walker, S.~G. (2019).
\newblock General {B}ayesian updating and the loss-likelihood bootstrap.
\newblock {\em Biometrika}, 106(2):465--478.

\bibitem[Madan et~al., 1998]{madan1998variance}
Madan, D.~B., Carr, P.~P., and Chang, E.~C. (1998).
\newblock The variance gamma process and option pricing.
\newblock {\em Review of Finance}, 2(1):79--105.

\bibitem[Madan and Seneta, 1990]{madan1990variance}
Madan, D.~B. and Seneta, E. (1990).
\newblock The variance gamma ({VG}) model for share market returns.
\newblock {\em Journal of {B}usiness}, pages 511--524.

\bibitem[Nickl and S\"{o}hl, 2019]{MR4013745}
Nickl, R. and S\"{o}hl, J. (2019).
\newblock Bernstein-von {M}ises theorems for statistical inverse problems {II}:
  compound {P}oisson processes.
\newblock {\em Electronic Journal of Statistics}, 13(2):3513--3571.

\bibitem[Sato, 2013]{MR3185174}
Sato, K. (2013).
\newblock {\em L\'{e}vy Processes and Infinitely Divisible Distributions},
  volume~68 of {\em Cambridge Studies in Advanced Mathematics}.
\newblock Cambridge University Press, Cambridge.
\newblock Translated from the 1990 Japanese original, Revised edition of the
  1999 English translation.

\bibitem[Shen and Ghosal, 2015]{MR3426318}
Shen, W. and Ghosal, S. (2015).
\newblock Adaptive {B}ayesian procedures using random series priors.
\newblock {\em Scandinavian Journal of Statistics}, 42(4):1194--1213.

\bibitem[Short, 2013]{short2013improved}
Short, M. (2013).
\newblock Improved inequalities for the {P}oisson and binomial distribution and
  upper tail quantile functions.
\newblock {\em International {S}cholarly {R}esearch {N}otices}, 2013.

\bibitem[Syring and Martin, 2019]{MR3949316}
Syring, N. and Martin, R. (2019).
\newblock Calibrating general posterior credible regions.
\newblock {\em Biometrika}, 106(2):479--486.

\bibitem[Syring and Martin, 2021]{syring2021gibbs}
Syring, N. and Martin, R. (2021).
\newblock Gibbs posterior concentration rates under sub-exponential type
  losses.
\newblock {\em arXiv preprint {\tt arXiv:2012.04505}}.

\bibitem[Taleb, 2007]{taleb2007black}
Taleb, N.~N. (2007).
\newblock {\em The Black Swan: The Impact of the Highly Improbable}, volume~2.
\newblock Random house.

\bibitem[Ueltzh\"{o}fer and Kl\"{u}ppelberg, 2011]{adaptive.levy.2011}
Ueltzh\"{o}fer, F. A.~J. and Kl\"{u}ppelberg, C. (2011).
\newblock An oracle inequality for penalised projection estimation of
  {L}\'{e}vy densities from high-frequency observations.
\newblock {\em Journal of Nonparametric Statistics}, 23(4):967--989.

\bibitem[Walker, 2013]{MR3082220}
Walker, S.~G. (2013).
\newblock Bayesian inference with misspecified models.
\newblock {\em Journal of Statistical Planning and Inference},
  143(10):1621--1633.

\bibitem[Wu and Martin, 2020]{gpc.compare}
Wu, P.-S. and Martin, R. (2020).
\newblock A comparison of learning rate selection methods in generalized
  {B}ayesian inference.
\newblock {\tt arXiv:2012.11349}.

\bibitem[Zhang, 2006a]{MR2291497}
Zhang, T. (2006a).
\newblock From {$\epsilon$}-entropy to {KL}-entropy: analysis of minimum
  information complexity density estimation.
\newblock {\em The Annals of Statistics}, 34(5):2180--2210.

\bibitem[Zhang, 2006b]{MR2241190}
Zhang, T. (2006b).
\newblock Information-theoretic upper and lower bounds for statistical
  estimation.
\newblock {\em IEEE Transactions on Information Theory}, 52(4):1307--1321.

\end{thebibliography}

\end{document}